\documentclass[11pt]{article}

\usepackage[colorlinks,citecolor=blue,urlcolor=black]{hyperref}
\usepackage{bm}
\usepackage{amssymb}
\usepackage{amsmath}
\usepackage{amsthm}
\usepackage{mathrsfs}
\usepackage{dsfont}
\usepackage{thmtools} 
\usepackage{thm-restate}
\usepackage{booktabs}
\usepackage{multirow}
\usepackage{adjustbox}
\usepackage{xcolor}
\usepackage{natbib}
\usepackage{url}
\usepackage[title]{appendix}
\usepackage[doublespacing]{setspace}
\usepackage[margin=1in]{geometry}
\usepackage{empheq}  %equations with left brace

\usepackage{apptools}
\AtAppendix{\counterwithin{lemma}{section}}

\usepackage{enumerate}
\usepackage{enumitem}
\usepackage[title]{appendix}

\def\E{\mathrm{E}}

\def\mathR{\mathbb{R}}

\DeclareMathOperator*{\argmax}{\arg\max}

\newtheorem{theorem}{Theorem}[section]

\newtheorem{lemma}{Lemma}
\newtheorem{proposition}[theorem]{Proposition}

\newtheorem{example}{Example}

\newcommand{\blind}{1}

\begin{document}

\def\spacingset#1{\renewcommand{\baselinestretch}%
{#1}\small\normalsize} %\spacingset{1}

\if1\blind
{
	\title{Extremal Dependence of Moving Average Processes Driven by Exponential-Tailed L\'evy Noise}
	\author{Zhongwei Zhang\footnotemark[1]\thanks{Research Center for Statistics, University of Geneva, 1205 Geneva, Switzerland}, David Bolin\footnotemark[2]\thanks{Statistics Program, CEMSE Division, King Abdullah University of Science and Technology (KAUST), Thuwal 23955-6900, Saudi Arabia}, Sebastian Engelke\footnotemark[1], Rapha\"el Huser\footnotemark[2]}
	\date{}
	\maketitle
} \fi

\if0\blind
{
	\bigskip
	\bigskip
	\bigskip
	\begin{center}
		{\LARGE\bf Title}
	\end{center}
	\medskip
} \fi
\medskip

\begin{abstract}
Moving average processes driven by exponential-tailed L\'evy noise are important extensions of their Gaussian counterparts in order to capture deviations from Gaussianity, more flexible dependence structures, and sample paths with jumps.
Popular examples include non-Gaussian Ornstein--Uhlenbeck processes and type G Mat\'ern stochastic partial differential equation random fields.
This paper is concerned with the open problem of determining their extremal dependence structure. We leverage the fact that such processes admit approximations on grids or triangulations that are used in practice for efficient simulations and inference. 
These approximations can be expressed as special cases of a class of linear transformations of independent, exponential-tailed random variables, that bridge asymptotic dependence and independence in a novel, tractable way. 
This result is of independent interest since models that can capture both extremal dependence regimes are scarce and the construction of such flexible models is an active area of research.
This new fundamental result allows us to show that the integral approximation of general moving average processes with exponential-tailed L\'evy noise is asymptotically independent when the mesh is fine enough. Under mild assumptions on the kernel function we also derive the limiting residual tail dependence function.
For the popular exponential-tailed Ornstein--Uhlenbeck process we prove that it is asymptotically independent, but with a different residual tail dependence function than its Gaussian counterpart.
%The computational advantage of the SPDE-based formulation of non-Gaussian processes is thus readily applicable to modeling spatial extremes.
Our results are illustrated through simulation studies.
\end{abstract}

%\spacingset{1.4} % DON'T change the spacing!
%\clearpage

%%%%%%%%%%%%%%%%%%%%%%%%%%%%%%%%%%%%%%%%%%%%%%%%%%%%%%%%%%%%%%%%%%%%%
\section{Introduction}
Moving average processes, also referred to as process convolutions, are popular and natural constructions for non-stationary and non-Gaussian processes that are widely applied in spatial statistics \citep{Higdon2002, Cressie2002, Rodrigues2010,VerHoef2010}.
They are defined as 
\begin{equation}\label{mov_ave_process}
u(\bm{s}) = \int_\mathcal{D} G(\|\bm{s}-\bm{t}\|) \mathcal{M}({\mathrm d}\bm{t}), %\quad \bm{s}\in\mathcal{D},
\end{equation}
where $\mathcal{D}$ is a Borel subset of $\mathR^d$ that possibly can depend on $\bm{s}$, $G: \mathR_+ \rightarrow \mathR_+$ is a measurable $\mathcal{M}$-integrable function, and $\mathcal{M}$ is a random measure~\citep{Kallenberg2017}.
A classical example in the temporal domain $d=1$ is the non-Gaussian Ornstein--Uhlenbeck (OU) process developed in \citet{BarndorffNielsen2001} to model stochastic volatility in financial econometrics; see also an application of this model to longitudinal data in \citet{Asar2020}.

Another related approach for constructing non-Gaussian processes is via a stochastic partial differential equation (SPDE) \citep{BarndorffNielsen2001, Bolin2014, BolinWallin2020}. The stationary solution (if it exists) to such an SPDE admits an integral representation of the form~\eqref{mov_ave_process}.
In the spatial domain $\mathR^2$, an important instance of these constructions are the type G Mat\'ern SPDE random fields \citep{Bolin2014, BolinWallin2020}, a non-Gaussian extension of the popular Gaussian Mat\'ern SPDE random fields \citep{Lindgren2011}.

\begin{figure}[!t]
    \centering
    \includegraphics[width=0.8\textwidth]{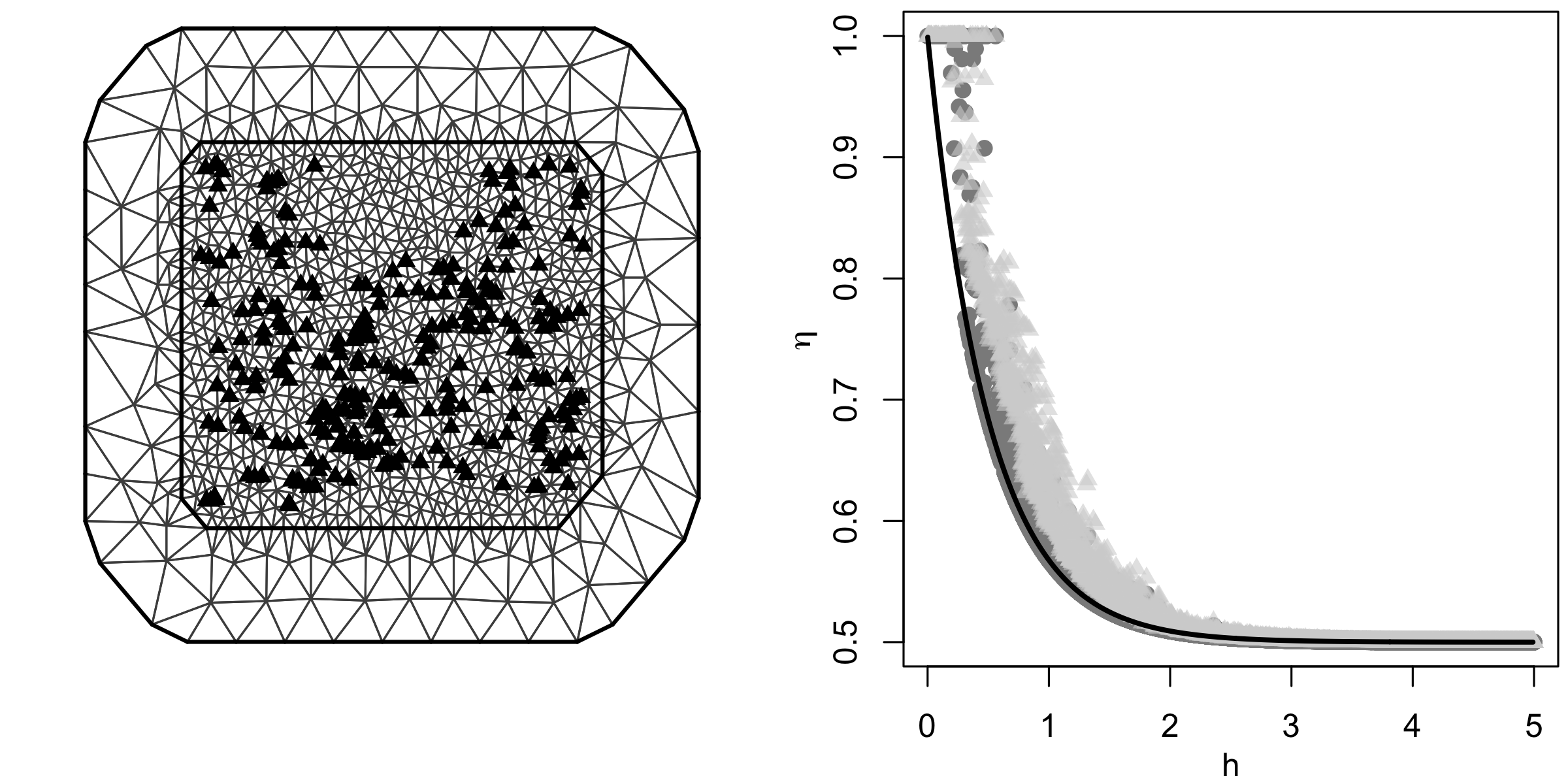}
    \caption{Observation locations and triangulation mesh used for Argo float data in \citet{BolinWallin2020} (left panel), and its induced residual tail dependence coefficients  $\eta$ for all pairs of locations for the finite element (light-grey triangle points) and integral approximations (dark-grey round points) of a type G Mat\'ern SPDE model, and their limiting function (black line).}
    \label{figure::Intro_argo}
\end{figure}

The key advantage of the SPDE-based process formulation is that one can approximate its solution by using the finite element method to obtain sparsity in the resulting precision matrix (or inverse of the dispersion matrix in the non-Gaussian case), thus achieving computationally efficient simulation and inference \citep{Lindgren2011, Bolin2014, BolinWallin2020}.
This finite element approximation has the form of a linear transformation
$
u_n(\bm{s}) = \sum_{i=1}^n a_i(\bm{s}) Y_i,
$
where the coefficients $a_i(\bm{s}) \geq 0$ are determined by the constructed basis functions on the triangulation with $n$ mesh nodes, and $Y_i$ are independent Gaussian or non-Gaussian random variables, depending on the random measure $\mathcal M$.
Non-Gaussian SPDE models are widely used in applications whenever the assumption of Gaussianity is too stringent. As an illustration, the left panel of Figure~\ref{figure::Intro_argo} shows the triangulation constructed in \citet{BolinWallin2020} for Argo float data, a data set that consists of measurements of seawater temperature and salinity in the global ocean and that has motivated the use of non-Gaussian processes due to its non-Gaussian features such as skewness and heavier tails than Gaussian distributions \citep{Kuusela2018}.

The dependence structure of SPDE and integral-based process constructions has mostly been studied in terms of their Pearson's correlation. This linear dependence measure does, however, not fully describe the process for non-Gaussian distributions. In fact, the correlation function of these processes is the same regardless of whether Gaussian or non-Gaussian noise is used.  
%, although it is a natural dependence measure in the world of Gaussian distributions.
%For instance, if one would like to answer questions such as (i) given the non-Gaussian OU process $u(t_1)$ is extreme, can $u(t_2)$ be also extreme at a later time point $t_2>t_1$? (ii) can the type-G Mat\'ern random field $u(\bm{s})$ be simultaneously extreme at two different locations $\bm{s}_1$ and $\bm{s}_2$? one has to understand the extremal dependence induced by these non-Gaussian processes.
Extremal (or tail) dependence describes the strength of dependence in the joint upper or lower tail of a multivariate distribution. It is crucial for risk assessment as it quantifies whether the largest realizations at different locations occur simultaneously.  
%, and was leveraged as a tool in \citet{Joe1993} to differentiate between different parametric copulas that lie between the independence copula and the comonotonicity copula.
This paper is concerned with the open problem of determining the extremal dependence structure of popular non-Gaussian processes and their linear approximations that are used in practice.
%As an illustration, consider the Argo measurements data application conducted in \citet{BolinWallin2020} for example, the constructed mesh is plotted on the left panel of Figure~\ref{figure::Intro_argo}.
%Here we are interested in the extremal dependence structure of the finite element approximation and its limiting form when the mesh size decreases to zero.

Let $(X_1,X_2)$ be a random vector with marginal distribution functions $F_{X_1}$ and $F_{X_2}$, respectively. A commonly used measure of extremal dependence is the (upper) tail dependence coefficient \citep{CHT1999}
\begin{equation}\label{chimeasure}
\chi = \lim_{q\uparrow 1} \chi(q) = \lim_{q\uparrow 1} \Pr(F_{X_1}(X_1) > q \mid F_{X_2}(X_2) > q),
\end{equation}
which satisfies $\chi\in [0,1]$ provided that this limit exists.
Since the extremal dependence in the lower tail can be obtained by negating the random vector, here we only focus on the upper tail.
We call $X_1$ and $X_2$ asymptotically dependent if $\chi>0$, and asymptotically independent otherwise.
%Asymptotic independence means that the joint tail decays faster than the marginal tail and it can be loosely interpreted as that $X_1$ and $X_2$ cannot be extreme simultaneously, whilst asymptotic dependence means that the joint tail decays at the same rate as the marginal tail and $X_1$ and $X_2$ can be simultaneously extreme.
In the latter case, the residual tail dependence coefficient $\eta$ \citep{Ledford1996} measures the second-order extremal dependence behavior and is defined as
\begin{equation}\label{etameasure}
\eta = \lim_{q\uparrow 1}\frac{\log(1-q)}{\log \Pr(F_{X_1}(X_1) > q, F_{X_2}(X_2) > q)},
\end{equation}
provided that this limit exists. 
For $\chi = 0$, the coefficient $1/\eta$ describes the rate of convergence of $\chi(q) \to 0$ as $q\to 1$. 
%The distinction between these two different extremal dependence regimes is important when evaluating the risk of joint extreme events since any misspecification could lead to either under- or overestimation of the risk.

%Therefore, to unify the framework, here we consider a general process
%\begin{equation}\label{exactModel}
%    X(\bm{s}) = \int_{\mathR^d} f(\bm{s}, \bm{t}) \Lambda({\mathrn d} \bm{t}), \quad \bm{s}\in\mathR^d,
%\end{equation}
%where $\Lambda$ is a L\'evy bases with exponential tails, i.e., for any bounded Borel set $A\in\mathR^d$, the distribution of $\Lambda(A)$ has an exponential tail, and $f: \mathR^d\times \mathR^d \rightarrow \mathR$ is a function satisfying certain conditions such that the stochastic integral is well-defined (conditions to be specified...).
The extremal dependence structure of continuously indexed processes of the form~\eqref{mov_ave_process} is challenging to derive due to the lack of analytical expressions of the induced multivariate distribution or density functions.
We therefore first consider a discretely indexed model of linear transformations
\begin{equation}\label{discretemodel}
    \left\{
     \begin{array}{cc}
       X_1 = a_{11}Y_1 + \cdots + a_{1n}Y_n, \\
       X_2 = a_{21}Y_1 + \cdots + a_{2n}Y_n,
     \end{array}
     \right.
\end{equation}
where $Y_i, i=1,\dots,n$, are independent and have exponential tails with the same index, and with coefficients $a_{ji} \geq 0,j=1,2, i=1,\dots,n$.   
%\footnote{I would delete "To avoid trivial cases, we assume that ${\max a_{1i} > 0}$, $\max a_{2i} > 0$, and $a_{1i}+a_{2i}>0$.", or this of importance here??}
This model is motivated by the finite element approximation of the Mat\'ern SPDE fields and integral approximation of general moving average processes.
For example, if we approximate~\eqref{mov_ave_process} by ${u_n(\bm{s})=\sum_{i=1}^n G(\|\bm{s}-\bm{t}_i\|)\mathcal{M}(D_i)}$ for a partition $\{D_i\}$ of $\mathcal{D}$ and ${\bm{t}_i\in D_i}$, and let ${X_j=u_n(\bm{s}_j)}$, then ${a_{ji}=G(\|\bm{s}_j-\bm{t}_i)}$ and ${Y_i=\mathcal{M}(D_i)}$.
The exponential tail of the noise variables appears naturally in the commonly used non-Gaussian processes.   
Under mild assumptions we show that $X_1$ and $X_2$ are asymptotically dependent if we have ${\argmax_{i\in\{1,\dots,n\}} a_{1i}=\argmax_{i\in\{1,\dots,n\}} a_{2i}}$, and they are asymptotically independent if ${\argmax_{i\in\{1,\dots,n\}} a_{1i} \cap \argmax_{i\in\{1,\dots,n\}} a_{2i} = \emptyset}$.
This novel result is of independent interest since it provides a tractable way to bridge asymptotic dependence and independence by simply varying the coefficients in (\ref{discretemodel}). It thus gives a partial answer to the second open problem raised in \citet{NoldeZhou2021} who ask how to build flexible parametric models that bridge both extremal dependence regimes.
%To the best of our knowledge, this is the first general result in the literature concerning the extremal dependence structure of linear transformations of independent random variables.

With this fundamental building block, we then study the extremal dependence structure of general moving average processes with exponential-tailed noise.
We show that under mild assumptions on the kernel function and the noise, the integral approximation of such processes is asymptotically independent when the mesh is fine enough, and we derive the limiting residual tail dependence function.
Moreover, we prove that the exponential-tailed non-Gaussian OU process \citep{BarndorffNielsen2001} is asymptotically independent, but with a different residual tail dependence function than its Gaussian counterpart.
%We further show that the integral approximation of the type G Mat\'ern SPDE model is asymptotically independent provided that the mesh is fine enough, and as the mesh size decreases to zero, the limiting residual tail dependence function is derived for certain values of the smoothness parameter $\alpha$. 
%We conjecture that asymptotic independence is preserved in the limiting process. 
As for the Argo float measurements data application, we consider the normal inverse Gaussian (NIG) Mat\'ern SPDE fields, which are a specific class of type G Mat\'ern SPDE fields.
%\footnote{"with SPDE parameters $\kappa=2, \alpha=3$": skip these details here but refer back to the picture later when we talk about results for SPDEs}
The right panel of Figure~\ref{figure::Intro_argo} shows the residual tail dependence coefficients of the finite element and integral approximations, and its theoretical limiting function as the triangulation mesh size goes to zero.
We conduct additional empirical studies in Section~\ref{TypeG_section} and Appendix~\ref{Appendix_simulation} to illustrate our results.

For moving average processes of the form~\eqref{mov_ave_process}, \citet{Opitz_unp2018} studied the extremal dependence structure when the kernel function $G$ is an indicator function, while \citet{Krupskii2022} studied the extremal dependence structure when $\mathcal{M}$ is a Cauchy random measure; for other cases, the extremal dependence structure remains unknown.
Our results are thus interesting since they motivate new constructions of non-Gaussian moving average processes and SPDE-based processes for extremes.

%Extreme value theory for processes of the form $X(\bm{s}) = \int_{\mathR^d} f(\bm{s}, \bm{t}) \mathcal{M}({\mathrm d} \bm{t})$ has been studied in the literature; see \citet{Rootzen1978}, \citet{Leadbetter1983}, \citet{Rootzen1986}, and more recently \citet{Stehr2021}.
%However, these works mainly focus on marginal extremes and they are extensions of the classical extremal types theorem \citep{FisherTippett1928,Gnedenko1943} from sequences of independent random variables to dependent sequences and then to continuously indexed processes.
%In contrast, our interest is in the extremal dependence structure induced by these processes, which, to our knowledge, has not been studied before.
%From the viewpoint of statistical modeling, our probabilistic results provide new insights into the effect of the mesh construction on the resulting extremal dependence of the type G Mat\'ern SPDE fields, which might provide some guidance on the mesh construction for practitioners when the interest is in modeling extremes.
%Moreover, we, for the first time, used a non-Gaussian SPDE model to describe the extremal dependence of a real dataset.

%The paper is structured as follows.
%Section 2 presents the necessary background knowledge on functions with exponential and regularly varying tails.
%The extremal dependence of the discretely indexed model, i.e., linear transformations of the form \eqref{discretemodel}, is discussed in Section 3, followed by that of the non-Gaussian OU processes in Section 4 and type G Mat\'ern SPDE fields in Section 5.
%Section 6 concludes with a discussion.

\section{Preliminaries on Exponential-Tailed Functions}\label{ExpoIntro}

A function $h: \mathR \rightarrow \mathR_+$ with $\lim_{x\rightarrow\infty} h(x) = 0$ is said to have an exponential tail with index $\beta\geq 0$, denoted as $h \in \mathscr{L}_{\beta}$, if 
\[
\lim_{x\rightarrow\infty}h(x-t)/{h(x)} = \exp(t\beta), \quad t\in\mathbb R.
\]
A univariate distribution function $F$ is said to have an exponential tail if its survival function $\bar{F}$ satisfies $\bar{F} \in \mathscr{L}_{\beta}$ for some $\beta\geq 0$. 
This is different from the conventional notation $F \in \mathscr{L}_{\beta}$, where the set $\mathscr{L}_{\beta}$ is restricted to be the family of all exponential-tailed distribution functions.
Our definition of $\mathscr{L}_{\beta}$ for general, nonnegative functions allows us to also cover exponential-tailed density functions.
If $\Bar{F} \in \mathscr{L}_{\beta}$ is such that we further have
$
\lim_{x\rightarrow\infty} {\overline{F*F}(x)}/{\bar{F}(x)} = M < \infty,
$
where $*$ denotes the convolution operator, then $F$ is called convolution tail equivalent, denoted by $\Bar{F}\in \mathscr{S}_\beta$.
Clearly, $\mathscr{S}_\beta \subset \mathscr{L}_{\beta}$.

A function $g: \mathR_+ \rightarrow \mathR_+$ is called regularly varying at $\infty$ with index $\rho\in\mathR$ if for any ${t>0}$,
$
{\lim_{x\rightarrow\infty} g(tx)/g(x) = t^\rho},
$
in which case we write $g\in {\text RV}_\rho$.
If $\rho=0$, $g$ is called slowly varying.
Exponential-tailed and regularly varying functions are intimately related to each other, since for $\beta \geq 0$, $h \in \mathscr{L}_{\beta}$ if and only if $h \circ \log \in {\text RV}_{-\beta}$. 
Using this relationship, we obtain a similar representation of exponential-tailed distribution functions to Karamata's representation of regularly varying functions (see Corollary 2.1 in \cite{Resnick2007}), namely for ${\Bar{F} \in \mathscr{L}_{\beta}}$, we have
\begin{equation}\label{expotailrepre}
{\bar F}(x) = a(x)\exp\Big\{-\int_0^x \beta(v) {\mathrm d}v  \Big\}, \quad x\in\mathR,
\end{equation}
where $a(x) \rightarrow a \in (0,\infty)$ and $\beta(x) \rightarrow \beta$ as $x\rightarrow\infty$. This will be repeatedly used in our proofs. 

Two prominent examples of exponential-tailed distributions are the generalized inverse Gaussian (GIG) distribution and the generalized hyperbolic (GH) distribution \citep{Barndorff1977, Prause1999, McNeil2005}.
A random variable $R$ is said to have a GIG distribution, denoted as $ R\sim \text{GIG}(\lambda,\tau,\psi)$, if its Lebesgue density is
\[
f_\text{GIG}(x) = \Big(\frac{\psi}{\tau}\Big)^{\lambda/2}\frac{x^{\lambda-1}}{2K_\lambda(\sqrt{\tau\psi})}\exp\Big\{-\frac{1}{2}\Big(\frac{\tau}{x}+\psi x\Big) \Big\}, \quad x>0,
\]
where $K_\lambda$ is the modified Bessel function of the second kind with index $\lambda$.
If a random variable $Z$ has the stochastic representation as a normal mean-variance mixture, i.e.,
\begin{equation*}
    Z = \mu + \gamma R + \sqrt{R}W, \quad \mu, \gamma\in\mathR,
\end{equation*}
where $W$ is a standard normal random variable and $R\sim \text{GIG}(\lambda,\tau,\psi)$, then $Z$ is said to have the GH distribution, denoted as $Z \sim \text{GH}(\lambda,\tau,\psi,\mu,\gamma)$, and its Lebesgue density is given in Appendix~\ref{Appendix_aux_results}.
%\[
%    f_{\mathrm GH}(x) = \frac{\psi^{\lambda/2} (\psi+\gamma^2)^{1/2 - \lambda}}{(2\pi)^{1/2}\tau^{\lambda/2}K_\lambda(\sqrt{\tau \psi})} \cdot \frac{K_{\lambda-1/2}\Big[\sqrt{\{\tau+(x-\mu)^2\}(\psi+\gamma^2)}\Big] e^{\gamma(x-\mu)}} {\Big[\sqrt{\{\tau+(x-\mu)^2\}(\psi+\gamma^2)} \Big]^{1/2 - \lambda}}, \quad x\in\mathR.
%\]
The admissible parameter values in both the GIG and GH distributions are ${\lambda<0, \tau>0, \psi\geq 0}$, or ${\lambda=0, \tau>0, \psi>0}$, or ${\lambda>0, \tau\geq 0, \psi>0}$.
The special cases $\psi=0$ and $\tau=0$ should be understood as limiting cases.
Another limiting case, when $R$ is degenerate and thus $Z$ follows a Gaussian distribution, is not considered in this paper as $Z$ does not have an exponential tail in this case.

It is easy to see that the GIG density function (thus, also its distribution function) has exponential tails.
The GH density and distribution functions also have exponential tails and the details are given in Appendix~\ref{Appendix_aux_results}.
%To see why the GH density function also has exponential tails, using the asymptotic expansion of the modified Bessel function of the second kind for large arguments \citep[Formula 9.7.2]{Abramowitz1972}, we obtain
%\begin{equation*}
%    f_{\mathrm GH}(x) \sim c_0 x^{\lambda-1}e^{-(\sqrt{\psi+\gamma^2}-\gamma)x}, \quad x\rightarrow \infty,
%\end{equation*}
%where $c_0>0$ is a constant, and $f(x)\sim g(x)$ as $x\rightarrow \infty$ means that $g$ is eventually non-zero and $f(x)/g(x) \rightarrow 1$ as $x\rightarrow \infty$.
%A simple application of L'H\^opital's rule gives that the GIG and GH distribution functions also have exponential tails.
Furthermore, the GIG and GH distributions are convolution tail equivalent if and only if $\lambda <0$; see \citet{Embrechts1982} for more details.
The GIG and GH distributions are infinitely divisible, and their associated L\'evy processes are widely used in finance; see \citet{Eberlein2001} and \citet{McNeil2005} for an overview.
The main examples in \citet{BarndorffNielsen2001} are the GIG L\'evy processes, while the non-Gaussian noise considered in \citet{Bolin2014}, \citet{WallinBolin2015}, and \citet{BolinWallin2020} is based on certain subclasses of the GH distribution.

%%%%%%%%%%%%%%%%%%%%%%%%%%%%%%%%%%%  Discrete Model  %%%%%%%%%%%%%%%%%%%%%%%%%%%%%%%%%%%%%%%%%
\section{Linear Transformations}\label{dissection}

\subsection{General Framework and Outline}
In this section, we focus on the extremal dependence of the random vector $\bm{X} = (X_1, X_2)$ defined in \eqref{discretemodel}
%, i.e., $X_1 = a_{11}Y_1 + a_{12}Y_2+\cdots+a_{1n}Y_n$, $X_2 = a_{21}Y_1 + a_{22}Y_2+\cdots+a_{2n}Y_n$, 
with $a_{ji}\geq 0$, $\max a_{1i} > 0$, $\max a_{2i} > 0$, and $a_{1i}+a_{2i}>0, j=1,2, i=1,\dots,n$.
We work under the following assumptions throughout the paper unless otherwise stated:
\begin{enumerate}[label=\textbf{A.\arabic*},ref=A.\arabic*]
    \item $Y_i, i=1,\dots,n$, are mutually independent with identical distribution function $F_{Y}$ such that $\Bar{F}_{Y} \in \mathscr{L}_\beta$ for some $\beta>0$;\label{assumptionA1}
    \item $F_{Y}$ is absolutely continuous with density $f_{Y}$. \label{assumptionA2}
\end{enumerate}
The assumption of identical distribution functions in condition \ref{assumptionA1} can be relaxed to different exponential-tailed distribution functions with the same index $\beta$, i.e., $Y_i\sim F_{Y_i}$ such that ${\Bar{F}_{Y_i} \in \mathscr{L}_\beta, \beta>0, i=1,\dots,n}$, but for the sake of simplicity we assume here that they have a common distribution $F_Y$.

Although the dependence structure in model (\ref{discretemodel}) is defined via simple linear transformations, surprisingly, nontrivial extremal dependence structures can be obtained.
Specifically, it turns out that the extremal dependence of $\bm{X}=(X_1,X_2)$ mainly depends on the largest coefficients among $a_{1i}, i=1,\dots,n$, and $a_{2i}, i=1,\dots,n$. 
We will show that the components of the vector $\bm{X}$ are asymptotically dependent if there is equality of the maximizing sets
\begin{align}\label{AD_case}
\argmax_{i\in\{1,\dots,n\}} a_{1i} = \argmax_{i\in\{1,\dots,n\}} a_{2i}. 
\end{align}
We discuss this case in Section~\ref{ADsubsection}.
On the other hand, there is asymptotic independence if these sets have an empty intersection 
\begin{align}\label{AI_case}
\argmax_{i\in\{1,\dots,n\}} a_{1i} \cap \argmax_{i\in\{1,\dots,n\}} a_{2i}=\emptyset, 
\end{align}
which is discussed in Section~\ref{AIsubsection}.
The delicate boundary case where the two maximizing sets ${\argmax_{i\in\{1,\dots,n\}} a_{1i}}$ 
and ${\argmax_{i\in\{1,\dots,n\}} a_{2i}}$ are not equal but have a non-empty intersection, is discussed in Section \ref{boundarycase}.

It is interesting to note that this behaviour is in sharp contrast to the case of heavy-tailed linear \citep[e.g.,][]{Gnecco2021} and max-linear \citep[e.g.,][]{Wang2011} models, where the random variables $Y_i$ have common survival function $\bar F_Y \in {\text RV}_{-\beta}$. 
In this more classical case, only asymptotic dependence (or complete independence) may arise and all coefficients $a_{ji}$ contribute to the corresponding tail dependence coefficient $\chi$. 

%%%%%%%%%%%%%%%%%%%%%%%%%%%%%%%%%%%  Asymptotic dependence
\subsection{Asymptotic Dependence}\label{ADsubsection}
Here we show that $\bm{X} = (X_1, X_2)$ defined via \eqref{discretemodel} is asymptotically dependent when~\eqref{AD_case} holds, and we give the explicit expression of its tail dependence coefficient.
All proofs are given in the Appendix.

\begin{proposition}\label{ADnY}
Let $Y_1,\dots,Y_n$ be a sequence of random variables satisfying the condition \ref{assumptionA1} with $\bar F_Y \in \mathscr{L}_{\beta}$ for some $\beta >0$, and let $\bm{X}=(X_1,X_2)$ be constructed as in \eqref{discretemodel}. 
If the set equality in~\eqref{AD_case} holds, then $X_1$ and $X_2$ are asymptotically dependent and the tail dependence coefficient of $\bm X$ can be expressed as
\begin{equation*}
    \chi = \E\left[\min\left\{\frac{\exp(\beta Z_1)}{M_{Z_1}(\beta)}, \frac{\exp(\beta Z_2)}{M_{Z_2}(\beta)} \right\}\right]
\end{equation*}
where $Z_1 = \sum_{i\notin I_{\text{max}}} a_{1i}Y_i/a_{1 \text{max}}$, $Z_2 = \sum_{i\notin I_{\text{max}}} a_{2i}Y_i/a_{2 \text{max}}$, and $M_{Z_j}$ is the moment generating function of $Z_j, j=1,2$. 
Here the index set $I_{\text{max}}$ is the common set of maximizers in~\eqref{AD_case} and $a_{j \text{max}}= \max_{i=1,\dots,n} a_{ji}$, $j=1,2$.
\end{proposition}

In this asymptotically dependent case where~\eqref{AD_case} holds, there is a link between our discrete model (\ref{discretemodel}) and the random scale constructions considered in \citet{Huser2019} and \citet{EOW2019}.
As shown in the proof of Proposition \ref{ADnY}, we can assume without loss of generality that $a_{1 \text{max}}=a_{2 \text{max}}=1$.
We can then rewrite model (\ref{discretemodel}) as $X_1=Z_c+Z_1$ and $X_2=Z_c+Z_2$, where $Z_c = \sum_{i\in I_{\text{max}}} Y_i$, and $Z_1$ and $Z_2$ are the same as in Proposition~\ref{ADnY}.
Further notice that $Z_c$ has an exponential tail if and only if $\exp(Z_c)$ has a regularly varying tail.
Hence, if one is interested in the extremal dependence structure of the random vector $(\exp(X_1), \exp(X_2))$, which is the same as that of $(X_1, X_2)$ since extremal dependence is invariant to monotonically increasing marginal transformations, then this problem falls into the setting of random scale constructions.
More precisely, we have  $(\exp(X_1), \exp(X_2)) = \exp(Z_c) (\exp(Z_1), \exp(Z_2))$ with $\exp(Z_c)$ being the shared random component.
An application of Proposition 1 in \citet{EOW2019} yields the asymptotic dependence of  $(\exp(X_1), \exp(X_2))$ and thus also $(X_1, X_2)$, and gives the same tail dependence coefficient as in Proposition \ref{ADnY}.
This alternative proof transforms the sum of exponential-tailed random variables into a product of regularly varying random variables and exploits the properties of regularly varying functions to derive the extremal dependence structure.
In comparison, our proof in Proposition \ref{ADnY} treats linear transformations of
a vector of independent exponential-tailed random variables in a direct manner and reveals many asymptotic properties of exponential-tailed distribution functions.
This might be of independent interest and these results could be useful in other contexts.

To further investigate the expression of the tail dependence coefficient $\chi$ in Proposition~\ref{ADnY}, we now consider the specific example where the $Y_i$ are GH distributed. 
Suppose ${Y_1, Y_2 \sim {\text{GH}}(\lambda, \tau, \psi, \mu, \gamma)}$ with $\psi>0, \gamma=0$, then Appendix~\ref{Appendix_aux_results} shows that their densities and survival functions have exponential tails with the same index, i.e., $f_{Y_i}, \bar F_{Y_i} \in \mathscr{L}_{\sqrt{\psi}}$.
Proposition~\ref{ADnY} yields the following result.

\begin{example}\label{ADtwoY}
Let $Y_1, Y_2$ be independent and have a common distribution ${\text{GH}}(\lambda, \tau, \psi, \mu, \gamma)$ with $\psi>0, \gamma=0$.
Let $X_1 = Y_1 + a_{12}Y_2$, $X_2 = Y_1 + a_{22}Y_2$ with $0\leq a_{12}, a_{22} < 1$.
Then $(X_1, X_2)$ is asymptotically dependent with tail dependence coefficient
\begin{equation*}
    \chi = \frac{\int_{(a_{22}-a_{12})y\leq c} \exp(a_{22}\sqrt{\psi}y)  F_{Y_1}({\mathrm d}y)}{M_{Y_1}(a_{22}\sqrt{\psi})} + 
    \frac{\int_{(a_{22}-a_{12})y> c} \exp(a_{12}\sqrt{\psi}y)  F_{Y_1}({\mathrm d}y)}{M_{Y_1}(a_{12}\sqrt{\psi})},
\end{equation*}
where $c=\{\log M_{Y_1}(a_{22}\sqrt{\psi}) - \log M_{Y_1}(a_{12}\sqrt{\psi}) \}/\sqrt{\psi}$, and $M_{Y_1}$ is the moment generating function of $Y_1$.
\end{example}
Clearly, the tail dependence coefficient $\chi$ in Example \ref{ADtwoY} is larger than zero.
Also note that when $a_{12}=a_{22}$, $X_1=X_2$ and thus $\chi=1$.
In fact, one can show that when $a_{12}\neq a_{22}$, $\chi$ is always strictly less than $1$; see Proposition \ref{Aux_result_1} in Appendix~\ref{Appendix_aux_results}.
To further investigate the properties of the tail dependence coefficient $\chi$ expressed in Example \ref{ADtwoY}, we now assume $0\leq a_{12} < a_{22} <1$ and examine its limit as $a_{22} \uparrow 1$.
This is interesting because in the following subsection we see that when $a_{12}<1<a_{22}$, which implies $\argmax_{i\in\{1,2\}} a_{1i} = 1$ and $\argmax_{i\in\{1,2\}} a_{2i}=2$, then $X_1$ and $X_2$ are asymptotically independent and necessarily $\chi=0$.
Hence, the investigation of $\lim_{a_{22}\uparrow 1} \chi$ answers the question of whether  $\chi$ smoothly transitions from asymptotic dependence to independence.
It turns out that this is true  only in some cases.
\begin{proposition}\label{ADtwoY_limit}
Let $\chi$ be as in Example \ref{ADtwoY} with ${0\leq a_{12} < a_{22} <1}$. Then
\begin{equation*}
    \lim_{a_{22}\uparrow 1} \chi = \begin{cases}
    \frac{\int_{-\infty}^{c^*} \exp(\sqrt{\psi}y)  F_{Y_1}({\mathrm d}y)}{M_{Y_1}(\sqrt{\psi})} + 
    \frac{\int^{\infty}_{c^*} \exp(a_{12}\sqrt{\psi}y)  F_{Y_1}({\mathrm d}y)}{M_{Y_1}(a_{12}\sqrt{\psi})}, & \text{ if } \lambda < 0,  \\
    0, & \text{ if } \lambda \geq 0,  
    \end{cases}
\end{equation*}
where $c^* = \lim_{a_{22}\uparrow 1} c/(a_{22}-a_{12}) = \frac{\log\{M_{Y_1}(\sqrt{\psi})\} - \log\{M_{Y_1}(a_{12}\sqrt{\psi}\}}{(1-a_{12})\sqrt{\psi}} < \infty$ when $\lambda < 0$.
\end{proposition}

\begin{figure}[!t]
    \centering
    \includegraphics[width=0.8\textwidth]{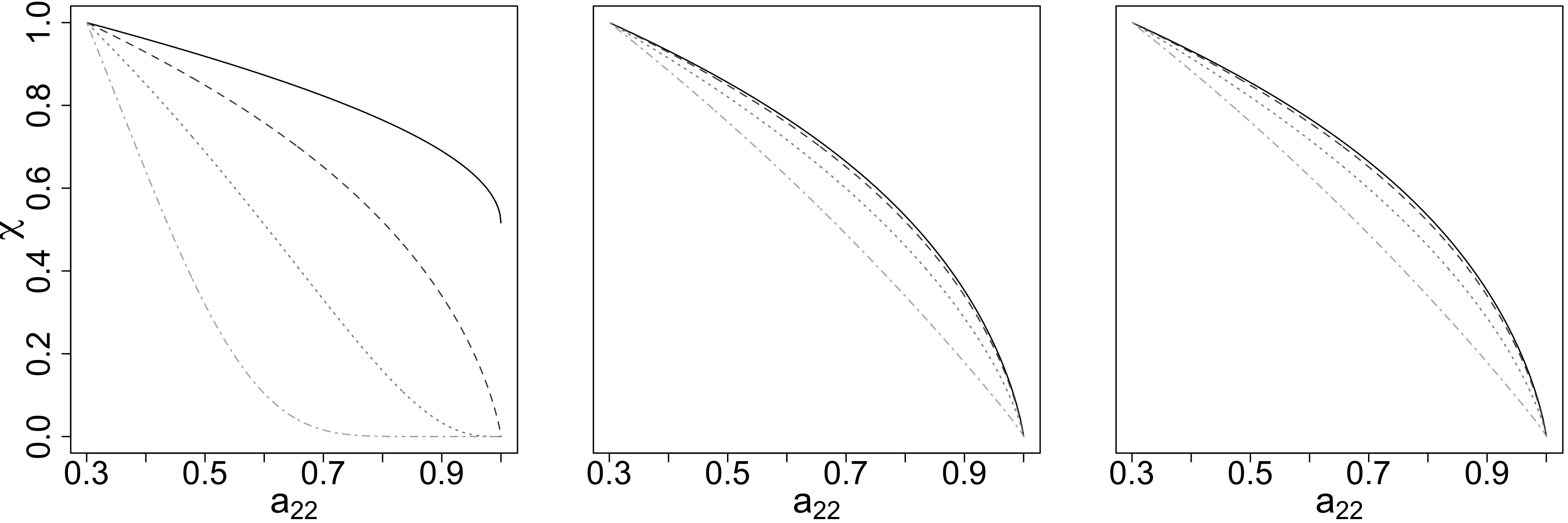}
    \caption{Values of $\chi$ as a function of $a_{22}$ for fixed $a_{12}=0.3$ and parameters $\tau=1, \psi=1, \lambda=-0.5, 1, 5, 30$ (left), $\lambda=1, \psi=1, \tau=0.5, 1, 5, 30$ (center), and $\lambda=1, \tau=1, \psi=0.5, 1, 5, 30$ (right) in Example~\ref{ADtwoY}, where lines in each plot are ordered from highest to lowest for increasing values of $\lambda, \tau$ and $\psi$, respectively.}
    \label{figure::limitchi_a22}
\end{figure}

The result in Proposition \ref{ADtwoY_limit} indicates that when $\lambda<0$, which implies that the distributions of $Y_1, Y_2$ are convolution tail equivalent \citep{Pakes2004}, there is a discontinuity in $\chi$ when $a_{22}$ tends to $1$ from below.
To illustrate how fast the $\chi$ coefficient in Example \ref{ADtwoY} tends to its limit, we set $a_{12}=0.3$, and plot $\chi$ against $a_{22}$ for various values of $\lambda,\tau,\psi$ in Figure \ref{figure::limitchi_a22}.
The results show that $\lambda$ is the most important parameter that regulates the decay rate of $\chi$ as $a_{22}$ increases, and the larger $\lambda$ the faster $\chi$ decays.
This indicates that the normal inverse Gaussian distribution, which is a subclass of the GH distribution when $\lambda$ is fixed to $-0.5$, might not be a good option to use for modeling extreme events in the presence of both asymptotic dependence and independence, since there is no smooth transition between the two. 
On the other hand, the variance gamma distribution, another subclass of the GH distribution with $\lambda > 0$, might be a better candidate to consider as there is a smooth transition.

%%%%%%%%%%%%%%%%%%%%%%%%%%%%%%%%%%%  Asymptotic independence
\subsection{Asymptotic Independence}\label{AIsubsection}
In this section we show that $\bm{X}$ defined in \eqref{discretemodel} is asymptotically independent when~\eqref{AI_case} holds.
We further give the explicit expression of the residual tail dependence coefficient $\eta$; recall its definition \eqref{etameasure} in the introduction.

The idea of this proof is to first expand the bivariate random vector $\bm{X}$ to $n$ dimensions, then use the geometric approach \citep{Nolde14,Nolde2022} to cast the computation of the residual tail dependence coefficient $\eta$ as a convex optimization problem, and finally solve this optimization problem.
We refer to the Supplementary Material (Section 1) for the proof and additional related results.

\begin{proposition}\label{AInY}
Let $Y_1,\dots, Y_n$ be a sequence of random variables satisfying~\ref{assumptionA1} and \ref{assumptionA2}.
Let $\bm{X}=(X_1,X_2)$ be constructed as in (\ref{discretemodel}).
Then the residual tail dependence coefficient is
\begin{align*}
   \eta = \Big[\min_{\substack{i,j=1,\dots,n \\ i\neq j}} \Big\{ \frac{|\Tilde{a}_{2i}-\Tilde{a}_{1i}|+|\Tilde{a}_{2j}-\Tilde{a}_{1j}|}{|\Tilde{a}_{2i}\Tilde{a}_{1j}-\Tilde{a}_{1i}\Tilde{a}_{2j}|}, \max(1/\Tilde{a}_{1i},1/\Tilde{a}_{2i}), \max(1/\Tilde{a}_{1j},1/\Tilde{a}_{2j}) \Big\} \Big]^{-1},
\end{align*}
where $\Tilde{a}_{ji}= a_{ji}/(\max_{r\in\{1,\dots,n\}}a_{jr}), j=1,2, i=1,\dots,n$, and the first term in the minimum is set to $+\infty$ whenever $\Tilde{a}_{2i}\Tilde{a}_{1j}-\Tilde{a}_{1i}\Tilde{a}_{2j}=0$.
Moreover, if~\eqref{AI_case} holds then $\eta<1$, and necessarily $X_1$ and $X_2$ are asymptotically independent.
If \eqref{AI_case} does not hold, then $\eta=1$.
\end{proposition}

We remark that when \eqref{AI_case} does not hold, more assumptions are needed to determine the extremal dependence regime of $\bm{X}$, i.e., whether the tail dependence coefficient $\chi>0$ or $\chi=0$; we refer to Section \ref{ADsubsection} for the asymptotic dependence case and Section \ref{boundarycase} for the boundary case.

We now consider the specific case that is similar to Example \ref{ADtwoY} and where a direct application of Proposition \ref{AInY} yields a simple form of $\eta$.

\begin{example}\label{AItwoY}
Let $Y_1, Y_2$ be independent and have a common distribution ${\text{GH}}(\lambda, \tau, \psi, \mu, \gamma)$ with ${\psi>0}$, ${\gamma=0}$.
Let $X_1 = Y_1 + a_{12}Y_2$, $X_2 = a_{21}Y_1 + Y_2$ with $0\leq a_{12}, a_{21} \leq 1$.
Then, the residual tail dependence coefficient of $\bm{X}=(X_1, X_2)$ is
\begin{equation*}
    \eta = \frac{1-a_{12}a_{21}}{2-a_{12}-a_{21}}.
\end{equation*}
Furthermore, if $a_{12}, a_{21}$ are both strictly less than $1$, then $\bm{X}$ is asymptotically independent.
\end{example}

An interesting observation in Proposition \ref{AInY} and Example \ref{AItwoY} is that the expression of $\eta$ only depends on the coefficients $a_{ji}$.
In other words, this means that as long as we do not change the coefficients $a_{ji}$ in the model (\ref{discretemodel}) and $Y_1,\dots,Y_n$ satisfy conditions \ref{assumptionA1} and \ref{assumptionA2}, then we always obtain the same residual tail dependence coefficient $\eta$, regardless of coefficient $\beta$ or the exact distribution of $Y_i$.
Another observation of Example \ref{AItwoY} is that if we fix $a_{12}$ (or $a_{21}$), then $\eta$ is monotonically increasing with respect to $a_{21}$ (or $a_{12}$).
%In the spatial context, this means that when we decrease the distance between two points, then $\eta$ increases to $1$.
Furthermore, the range of $\eta$ is $[1/2, 1]$, where $1/2$ is achieved when $a_{12}=a_{21}=0$, and $\eta=1$ when $a_{12}=1$ or $a_{21}=1$.

We now consider another example that clearly distinguishes the extremal dependence measures from dependence measures for the bulk of the distribution.

\begin{example}\label{AIthreeY}
Let $Y_1, \dots, Y_n$ be independent and identically distributed with common distribution function $F_Y$ and $\Bar{F}_{Y} \in \mathscr{L}_{\beta}, \beta>0$.
Let ${X_1 = Y_1}$, ${X_2 = a_{21}Y_1 + Y_2 + a_{23}Y_3+\cdots+a_{2n}Y_n}$ with ${0\leq a_{21}, a_{23},\dots, a_{2n} < 1}$.
Then by Proposition~\ref{AInY}, $\bm{X}=(X_1, X_2)$ is asymptotically independent with residual tail dependence coefficient $\eta = 1/(2-a_{21})$.
\end{example}

%Notice that when $a_{11}=a_{21}=a_0 < 1$, the residual tail dependence coefficient reduces to $\eta = \max(a_0, 1/(2-a_{22}))$.
%If we denote $R=\exp(a_0 Y_1)$, $W_1=\exp(Y_2)$, and $W_2=\exp(a_{22}Y_2+Y_3)$, then from Example \ref{AItwoY} we know that $W_1$ and $W_2$ are asymptotically independent and the corresponding residual tail dependence coefficient is $1/(2-a_{22})$.
%By Proposition 5 in \citet{EOW2019}, one obtains the residual tail dependence coefficient of the random scale construction $(\exp(X_1), \exp(X_2))=R (W_1, W_2)$ as $\max(a_0, 1/(2-a_{22}))$, which is identical to that of $(X_1, X_2)$.
We observe that the expression of $\eta$ does not depend on $n$, i.e., the number of independent terms in the construction of $X_2$.
Loosely speaking, when more independent terms ${a_{2i}Y_i}$ with ${0< a_{2i}< 1}$ are added to $X_2$, these added terms do not contribute much to the extremal dependence between $X_1$ and $X_2$ and the residual tail dependence coefficient remains the same.
This is in clear contrast with the classical Pearson's correlation dependence measure
\[ \text{Corr}(X_1,X_2)={a_{21}}\Big({a_{21}^2+1+\sum_{i=3}^n a_{2i}^2}\Big)^{-1/2}, \]
which decreases monotonically as $n$ increases.

%%%%%%%%%%%%%%%%%%%%%%%%%%%%%%%%%%%  Remarks
\subsection{Boundary Case and Further Remarks}\label{boundarycase}
In Sections \ref{ADsubsection} and \ref{AIsubsection} we have shown that the discrete model $\bm{X}$ in \eqref{discretemodel} is asymptotically dependent or independent when~\eqref{AD_case} or~\eqref{AI_case} hold, respectively.
The intuition for this phenomenon is that the extremal dependence of $(X_1,X_2)$ is determined by whether the main contribution for the tails of $X_1$ and $X_2$ comes from the same components $Y_i$. 
If this is the case, then the stronger form of extremal dependence, namely asymptotic dependence, is achieved; otherwise, asymptotic independence is obtained. In the sequel, let $I_j = \argmax_{i\in\{1,\dots,n\}} a_{ji}$, $j=1,2$.

Clearly, there is a boundary case where neither~\eqref{AD_case} nor~\eqref{AI_case} hold, that is, we have $I_1 \neq I_2$ and  
\[ I_{\text{max}} := I_1  \cap I_2\neq \emptyset.\]
To avoid complications, we here only consider the case where $I_1 \not\subset I_2 $ and $I_2 \not\subset I_1$.
With the notation of Proposition~\ref{ADnY} and the paragraph thereafter, without losing generality, we assume that $a_{1 \text{max}}= a_{2 \text{max}} = 1$, and recall that we can write $(X_1,X_2)=Z_c+(Z_1,Z_2)$, or equivalently, $(\exp(X_1), \exp(X_2)) = \exp(Z_c) (\exp(Z_1), \exp(Z_2))$.
Theorem 3 in \citet{Embrechts1980} now implies that the survival functions of $\exp(Z_c)$, $\exp(Z_1)$, and $\exp(Z_2)$ are all regularly varying with the same index $-\beta$.
Hence, Proposition 6 in \citet{EOW2019} can be applied to determine the extremal dependence structure of $(X_1,X_2)$.
More precisely, in this case, the residual tail dependence coefficient of $(X_1,X_2)$ is $\eta=1$ (this can also be obtained using Proposition~\ref{AInY}), but more assumptions are needed to determine whether $\chi>0$ or $\chi=0$.
We leave the other boundary case $I_\text{max} \neq \emptyset$, $I_1 \subset I_2$, or \emph{vice versa}, for future research.

Here we have studied the extremal dependence structure of linear transformations, or sums of exponential-tailed random vectors.
The results are directly applicable to products of regularly varying random vectors.
Indeed, assume that $\bar{Y}_1, \dots, \bar{Y}_n$ are independent copies of a positive random variable $Y$ with absolutely continuous distribution function $F_{\bar{Y}}$, and $\bar{F}_{\bar{Y}} \in {\text RV}_{-\beta}$ with $\beta>0$.
Let $a_{ji}\geq 0$ for $j=1,2, i=1,\dots,n$.
Then the product model 
\begin{equation*}
 \left\{
  \begin{array}{cc}
    \bar{X}_1=\bar{Y}_1^{a_{11}}\cdots \bar{Y}_n^{a_{1n}}, \\
    \bar{X}_2=\bar{Y}_1^{a_{21}}\cdots \bar{Y}_n^{a_{2n}},
  \end{array}
  \right.
\end{equation*}
 has the same extremal dependence structure as the sum model \eqref{discretemodel} with $Y_i=\log(\bar{Y_i})$, and Propositions \ref{ADnY} and \ref{AInY} give the extremal dependence coefficients of $(\bar{X}_1,\bar{X}_2)$. 

We here restrict our attention in~\eqref{discretemodel} to the case with non-negative coefficients $a_{ji}$, because both the integral approximation of moving average processes and the finite element approximation of the type G SPDE Mat\'ern fields have non-negative coefficients; see Section~\ref{mov_ave_section} for more details.
Consequently, the residual tail dependence coefficient $\eta$ in Proposition~\ref{AInY} can be shown to have the lower bound $1/2$, which implies that only positive extremal association can be achieved.
However, if the interest lies in capturing negative extremal association, i.e., $\eta<1/2$, then one can achieve this by considering negative coefficients $a_{ji}$ and assuming that both the left and right tails of the distribution of $Y_i$ have exponential decays.
We leave this for future research.

%%%%%%%%%%%%%%%%%%%%%%%%%%%%%%%%%%%%%%%% Continuous Model %%%%%%%%%%%%%%%%%%%%%%%%%%%%%%%%%%%%%%%%
\section{Moving Average Processes}\label{mov_ave_section}
\subsection{Main Results}
We now focus on the extremal dependence structure of moving average processes~\eqref{mov_ave_process}.
Note that here $u(\bm{s})$ refers to the usual definition of stochastic integrals, where in the first step one defines integrals of simple functions, and then a non-random integrand function $G$ is called $\mathcal{M}$-integrable if there exists a sequence of simple functions that converges pointwise to $G$, such that the limit in probability of the resulting integrals of simple functions exists, and this limit is defined to be the stochastic integral with respect to $G$ \citep{Rajput1989}.
A useful characterization of $\mathcal{M}$-integrable functions is given in \citet[Theorem 2.7]{Rajput1989}.
Throughout the paper we assume that $u(\bm{s})$ is well defined, i.e., $G$ is $\mathcal{M}$-integrable.

We first consider the more commonly used case where the domain of integration $\mathcal{D}$ is fixed and does not depend on $\bm{s}$.
To consider a framework that is applicable to general moving average processes which are not necessarily SPDEs, we only assume that: 
\begin{enumerate}[label=\textbf{B.\arabic*},ref=B.\arabic*]
    \item the function $G(h)$ is non-negative, continuous and strictly decreasing as $h\rightarrow\infty$; \label{assumptionB1}
    \item for any bounded Borel set $B\subset \mathcal D$, the random variable $\mathcal{M}(B)$ has an absolutely continuous distribution function $F_{\mathcal{M}(B)}$ with exponential tail $\Bar{F}_{\mathcal{M}(B)} \in \mathscr{L}_\beta$, $\beta>0$. \label{assumptionB2}
\end{enumerate}
We note that for the important class of type G Mat\'ern SPDE random fields on $\mathR^d$, the function $G$ is non-negative, absolutely continuous, and monotonically decreasing (see Proposition~\ref{Aux_result_3} in Appendix~\ref{Appendix_aux_results}).
Assumption~\ref{assumptionB2} implies that $\mathcal{M}(B)$ has an exponential-tailed distribution with the same index for every bounded Borel set $B$.
This seemingly restrictive assumption is, in fact, a natural result of the convolution-closure property of exponential tails \citep[Theorem 3]{Embrechts1980}, and Assumption~\ref{assumptionB2} is satisfied for instance for the NIG Mat\'ern SPDE fields and the variance Gamma Mat\'ern SPDE fields considered in \citet{Bolin2014}, \citet{WallinBolin2015} and \citet{BolinWallin2020}.

We now assume that $d=2$, which is the most important case for spatial applications; the cases $d=1$ and $d\geq 3$ can be treated analogously.
Let $\{D_n: D_n\subset \mathcal{D}\}$ be an increasing sequence of bounded Borel sets in $\mathR^2$, and $D_n^1, D_n^2, \dots, D_n^J \subset D_n$ be a partition of $D_n$ obtained by triangulating $D_n$ with $m_n$ mesh nodes $M_n=\{\bm{a}_n^i,i=1,\dots,m_n\}$.
Then we obtain an approximation of $u(\bm{s})$ as
\begin{equation}\label{mov_ave_approx}
    u_n(\bm{s}) = \int_{\mathcal{D}} \sum_{i=1}^J G(\|\bm{s}-\bm{d}_i\|) \mathrm{1}_{D_n^i}(\bm{t}) \mathcal{M}({\mathrm d} \bm{t}) = \sum_{i=1}^J G(\|\bm{s}-\bm{d}_i\|)\mathcal{M}(D_n^i),
\end{equation}
where $\bm{d}_i\in D_n^i$, and $\mathrm{1}_{D_n^i}$ is the indicator function.
We say that the sequence of points $M_n$ is dense in $\mathcal{D}$ if for any point $\bm{s}\in\mathcal{D}$, there is a sequence ${\{\Bar{\bm{a}}_n\}, \Bar{\bm{a}}_n\in M_n}$ such that ${\lim_{n\to\infty} \|\Bar{\bm{a}}_n - \bm{s}\|=0}$.
Assume that $M_n$ is dense.
Clearly, for any fixed $\bm{s}\in\mathcal{D}$, function ${G_n^J(\|\bm{s}-\bm{t}\|)=\sum_{i=1}^J G(\|\bm{s}-\bm{d}_i\|) \mathrm{1}_{D_n^i}(\bm{t})}$ converges to the function $G(\|\bm{s}-\bm{t}\|)$ pointwise as $m_n \rightarrow\infty$.
Hence, by the definition of stochastic integrals, we know that $u_n(\bm{s})$ converges to $u(\bm{s})$ in probability.
It follows that ${(u_n(\bm{s}_1), u_n(\bm{s}_2))}$, ${\bm{s}_1\neq \bm{s}_2}$, converges in probability to ${(u(\bm{s}_1),u(\bm{s}_2))}$.
%Indeed, a constructive definition of stochastic integrals is to approximate the integrand function by simple functions and then define the integral as the limit in probability; see Section 3.4 in \citet{Samorodnitsky1994} for such a definition of stable integrals.

Furthermore, by Assumption~\ref{assumptionB2}, ${\Bar{F}_{\mathcal{M}(D_n^i)} \in \mathscr{L}_\beta}$.
Note also that $\mathcal{M}(D_n^i), i=1,\dots,J$, are independent since $\{D_n^i, i=1,\dots,J\}$ forms a partition of $D_n$.
Hence, to understand the extremal dependence structure of ${(u_n(\bm{s}_1), u_n(\bm{s}_2))}$, it is sufficient to focus on the coefficients ${G(\|\bm{s}_1 - \bm{d}_i\|)}$ and ${G(\|\bm{s}_2 - \bm{d}_i\|)}$.
Since $\bm{s}_1\neq \bm{s}_2$, we know that if the mesh is fine enough, $\bm{s}_1$ and $\bm{s}_2$ will fall into different triangles $D_n^{i_1}$ and $D_n^{i_2}$.
Consequently, we have 
$$\argmax_i G(\|\bm{s}_1 - \bm{d}_i\|) = i_1 \neq \argmax_i G(\|\bm{s}_2 - \bm{d}_i\|) = i_2.$$
Proposition \ref{AInY} then yields the asymptotic independence of $(u_n(\bm{s}_1), u_n(\bm{s}_2))$ and allows the computation of its residual tail dependence coefficient. We focus on the limit of the residual tail dependence coefficient of the approximation model $(u_n(\bm{s}_1), u_n(\bm{s}_2))$ as $m_n\rightarrow\infty$.

\begin{theorem}\label{TypeG_limit_eta_prop}
Let $u_n(\bm{s})$ be the integral approximation of the moving average process~\eqref{mov_ave_process}, and let $\eta_n(h)$ be the residual tail dependence coefficient of $(u_n(\bm{s}_1), u_n(\bm{s}_2))$ for $h = \|\bm{s}_1 - \bm{s}_2\|$.
If Assumptions~\ref{assumptionB1} and~\ref{assumptionB2} are satisfied and the function $G$ is convex, then if $m_n\rightarrow\infty$ and the sequence of mesh nodes $M_n$ is dense in $\mathcal{D}$, the limit of $\eta_n$ as $n\to\infty$ is
\[
 \eta(h) = 
    \frac{1}{2}+\frac{G(h)}{2G(0)}.
\]
\end{theorem}

The idea of our proof is to use Proposition~\ref{AInY} to obtain an explicit formula for $\eta_n(h)$ as a minimum (or maximum) over a finite number of terms, and with this, to show that for convex $G$, the optimum can be achieved as the mesh size tends to zero.
We remark that when $G(0)=\infty$, $\eta(h)$ reduces to a constant $1/2$.
When the function $G$ is non-convex, we conjecture that the limiting residual tail dependence function is 
\[
\eta(h) = \max\Big\{\frac{1}{2}+\frac{G(h)}{2G(0)}, \frac{G(h/2)}{G(0)} \Big\}.
\]
Our simulation studies in Section~\ref{TypeG_section} seem to support our conjecture, but a rigorous proof would have to use a different technique than the proof of Theorem~\ref{TypeG_limit_eta_prop}, which relies heavily on the convexity of $G$.

We now consider another important case where the integration domain $\mathcal{D}$ depends on $s$, namely when $d=1$ and 
\begin{equation}\label{mov_OU_process}
{u(s)=\int_{-\infty}^s G(s-t) \mathcal{M}(\mathrm{d}t)}, \quad s \in (-\infty,T].
\end{equation}
This is an interesting case because it is used in practice \citep{BarndorffNielsen2001, VerHoef2010}, and this one-sided integral turns out to yield different residual tail dependence functions.

Let $-n< s_1 < s_2 \leq T$ and ${-n=t_0 < \cdots< t_{n_1} \leq s_1 < \cdots \leq t_{n_2} \leq s_2 < \cdots < t_{m_n} = T}$ be an arbitrary partition  of $[-n,T]$.
Then the approximation~\eqref{mov_ave_approx} becomes
\begin{equation}\label{mov_OU_approx}
  \begin{aligned}
    u_n(s_1) = \sum_{i=0}^{n_1-1} G(s_1 - t_i)\mathcal{M}([t_i, t_{i+1})), \quad 
    u_n(s_2) = \sum_{i=0}^{n_2-1} G(s_2 - t_i)\mathcal{M}([t_i, t_{i+1})). 
  \end{aligned}
\end{equation}
Clearly when $m_n\rightarrow\infty$ such that the sequence of partition points $M_n=\{t_i, i=0,\dots,m_n\}$ is dense in $(-\infty,T]$, ${(u_n(s_1), u_n(s_2))}$ converges in probability to ${(u(s_1), u(s_2))}$.
Note that when the mesh is coarse and there is no partition point in $(s_1,s_2)$, i.e., $t_{n_1}=t_{n_2}$, then the largest coefficients in the expressions of $u_n(s_1)$ and $u_n(s_2)$ are $G(s_1-t_{n_1-1})$ and $G(s_2-t_{n_1-1})$ respectively, which correspond to the same variable $\mathcal{M}([t_{n_1-1}, t_{n_1}))$.
Hence, Proposition~\ref{ADnY} gives the asymptotic dependence of ${(u_n(t_1), u_n(t_2))}$.
Otherwise, when there is at least one partition point in $(s_1,s_2)$, ${(u_n(t_1), u_n(t_2))}$ is asymptotically independent, and in the following, we derive its limiting residual tail dependence function as the mesh size tends to zero.

\begin{theorem}\label{OU_limit_eta}
Let $u_n(s)$ be the integral approximation~\eqref{mov_OU_approx} of the one-side integral~\eqref{mov_OU_process}, and let $\eta_n(h)$ be the residual tail dependence coefficient of $(u_n(s_1), u_n(s_2))$ for $h=|s_1 -s_2|$.
If Assumptions~\ref{assumptionB1} and~\ref{assumptionB2} are satisfied, then if $m_n\rightarrow\infty$ and the sequence of mesh nodes $M_n$ is dense in $\mathcal{D}$, the limit of $\eta_n$ as $n\to\infty$ is
\[
\eta(h) = \frac{1}{2-G(h)/G(0)}.
\]
\end{theorem}

\subsection{Non-Gaussian Ornstein-Uhlenbeck (OU) Process}\label{OUsection}

As the first application of the general results of the previous section, here we study the extremal dependence structure of non-Gaussian OU processes $u(t)$ \citep{BarndorffNielsen2001}, defined as the stationary solution to the stochastic differential equation (SDE)
\begin{equation}\label{OUmodel}
    {\mathrm d}u(t) = -a u(t){\mathrm d}t + {\mathrm d}z(a t), \quad a>0, t\in\mathR,
\end{equation}
where $z=\{z(t): t\in\mathR \}$ is a L\'evy process satisfying $\E[\log\{1+|z(1)|\}] < \infty$ to guarantee the existence of such a stationary solution.
Although the background driving L\'evy process $z$ can be chosen arbitrarily, examples considered in \citet{BarndorffNielsen2001} all have exponential tails for the density of the L\'evy measure of $z(1)$.

Non-Gaussian OU processes are generalizations of classical OU processes by replacing the Brownian motion $z(t)$ in the SDE~\eqref{OUmodel} by general L\'evy processes. 
The existence of such processes is established based on the notion of self-decomposability and the stochastic integral representation of self-decomposable random variables \citep{JurekVervaat1983}.
More precisely, let $V$ be a self-decomposable random variable (namely for every $\alpha\in(0,1)$, we have the decomposition $V=\alpha V+\Tilde{V}_\alpha$, where $\Tilde{V}_\alpha$ is a random variable independent of $V$), then there exists a L\'evy process $z(t)$ and a stationary stochastic process $u(s)$ such that \eqref{OUmodel} holds for all $a>0$ and $u(s)$ has the same distribution as $V$ for all $s \geq 0$.
Conversely, if $z$ is a L\'evy process and $u(t)$ is a stationary stochastic process such that $u(s)$ satisfying \eqref{OUmodel} for all $a>0$, then the marginal distribution of $u(s)$ is self-decomposable.
We refer to \citet{JurekVervaat1983} and \citet{BarndorffNielsen2001} for more details.

Importantly, the solution to the SDE~\eqref{OUmodel} can be represented as
\begin{equation}\label{OUrepresentation}
u(s) = \int_{-\infty}^s e^{-a(s-t)}{\mathrm d}z(at) = e^{-a s}u(0) + \int_0^s e^{-a(s-t)}{\mathrm d}z(at), \quad s\geq 0,
\end{equation}
where $u(0)$ is independent of $\int_0^s e^{-a(s-t)}{\mathrm d}z(at)$.
Using this representation, we can show that if the stationary solution $u(s)$ has an absolutely continuous and exponential-tailed distribution $F$ with $\bar{F} \in \mathscr{L}_\beta, \beta>0$, then the process $u(s)$ is asymptotically independent.

\begin{theorem}\label{OUprocess_AI}
    Let $u(s)$ be a non-Gaussian OU process defined as a stationary solution to \eqref{OUmodel} and $V$ has the stationary self-decomposable distribution.
    If the distribution function of $V$ is absolutely continuous with exponential tail $\bar{F}_V \in \mathscr{L}_\beta$, $\beta > 0$, then $(u(s_1), u(s_2))$ is asymptotically independent for $s_1\neq s_2$. The corresponding residual tail dependence coefficient is
    \[
     \eta = \frac{1}{2-e^{-a|s_2 - s_1|}}.
    \]
\end{theorem}

Note that when the L\'evy process $z$ in~\eqref{OUmodel} is a Brownian motion, we obtain the classical Gaussian OU process with the correlation between $u(s_1)$ and $u(s_2)$ as $e^{-a|s_2-s_1|}$ and residual tail dependence coefficient as $(1+e^{-a|s_2-s_1|})/2$.
Theorem \ref{OUprocess_AI} shows that when the marginal distribution of the non-Gaussian OU process $u(s)$ has an exponential tail, its induced extremal dependence structure is indeed different from the classical Gaussian OU process, although its correlation function remains the same.
More precisely, using a first-order Taylor expansion, one can observe that as $|s_2-s_1|$ tends to zero, the residual tail dependence coefficient $\eta$ of the Gaussian OU process increases to $1$ at a linear rate $a/2$, whilst $\eta$ of the non-Gaussian OU process increases to $1$ at a linear rate $a$.

Instead of specifying the self-decomposable distribution function of $u(s)$, one can alternatively define the non-Gaussian OU process by specifying the background driving L\'evy process $z$, which reduces to specifying the distribution of $z(1)$.
Although the relationship between the tail of $z(1)$ and that of $u(t)$ is unclear, $u(t)$ and $z(1)$ are closely linked by their L\'evy measures $U_{u}$ and $U_{z(1)}$ through the equation $U_{u}([x,\infty)) = \int_1^\infty s^{-1} U_{z(1)}([sx,\infty)) {\mathrm d}s$ \citep[Theorem 2.2]{BarndorffNielsen1998}.
%Although it is unclear whether an exponential tail of $z(1)$ implies an exponential tail for the marginal distribution of $u(t)$, neither is the inverse relation clear, $u(t)$ and $z(1)$ are closely linked by their L\'evy measures $U_{u}$ and $U_{z(1)}$ through the equation $U_{u}([x,\infty)) = \int_1^\infty s^{-1} U_{z(1)}([sx,\infty)) {\mathrm d}s$ \citep{BarndorffNielsen1998}.
Using this link, we show in the following that under mild assumptions, a convolution equivalent tail of $z(1)$ implies a convolution equivalent tail for $u(t)$.

\begin{proposition}\label{OUexpotail}
    Denote the distribution function of $z(1)$ and $u(t)$ by $F_{z(1)}$ and $F_{u}$, respectively.
    Suppose that the function $U_{z(1)}([x,\infty))$ is continuous in $x$ on $(0,\infty)$.
    If ${\Bar{F}_{z(1)}\in\mathscr{S}_\beta}$, ${\beta\geq 0}$, then ${\Bar{F}_{u}\in\mathscr{S}_\beta}$. Moreover, ${\Bar{F}_{z(t)}\in \mathscr{S}_\beta}$, where $F_{z(t)}$ denotes the distribution function of $z(t)$.
    Hence, all increments of the L\'evy process $z(t)$ have 
    exponential-tailed distribution functions with the same index $\beta$.
\end{proposition}

One example is given by $U_{z(1)}([x,\infty)) = c x^{-\epsilon}\exp(-\beta x)$, where $\epsilon>1$ and $c,\beta>0$; see Example 2.4.1 in \citet{BarndorffNielsen2001}.
The restriction $\epsilon>1$ implies that the normalized L\'evy measure of $z(1)$ is convolution tail equivalent
%, i.e., $\frac{\mathds{1}_{x>1} U_{z(1)}([x,\infty))}{U_{z(1)}([1,\infty))} \in\mathscr{S}_\beta$.
and hence $\Bar{F}_{z(1)} \in\mathscr{S}_\beta$ \citep[][Theorem B]{Shimura2005}.

Now we assume that $z(1)$ has a convolution tail equivalent distribution $F_{z(1)}$ with ${\Bar{F}_{z(1)}\in \mathscr{S}_\beta}$, ${\beta> 0}$, the function $U_{z(1)}([x,\infty))$ is continuous in $x$ on $(0,\infty)$, and the distribution function of $z(t)$ is absolutely continuous for all $t$.
That is, Assumption~\ref{assumptionB2} is satisfied.
Then we are ready to link our results in Theorem~\ref{OU_limit_eta} and Theorem~\ref{OUprocess_AI}.

\begin{figure}[!t]
    \centering
    \includegraphics[width=0.4\textwidth]{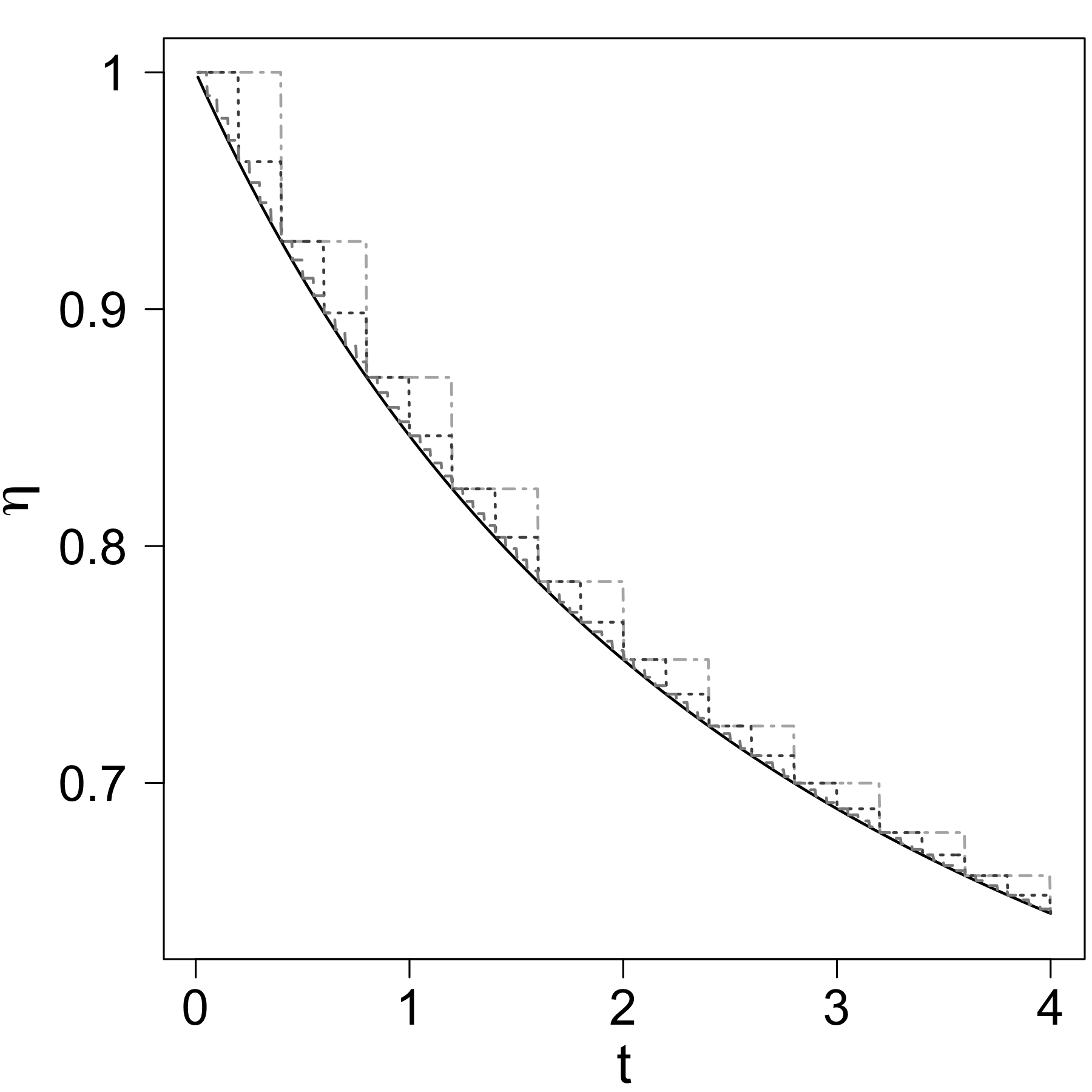}
    \caption{The residual tail dependence coefficient $\eta$ of the true non-Gaussian OU process (black line), and its integral approximations for three different, equidistant partitions of $[0,T]$ with mesh length $\Delta \in \{0.4, 0.2, 0.05\}$, where the lines are ordered from highest to lowest for decreasing $\Delta$.}
    \label{figure::limiteta_OU}
\end{figure}

Clearly, for the OU processes, the integrand function in~\eqref{OUrepresentation} satisfies Assumption~\ref{assumptionB1}.
Hence, Theorem~\ref{OU_limit_eta} implies that the limiting residual tail dependence function of the approximation model of the form~\eqref{mov_OU_approx} is $\eta(h)=1/(2-e^{-ah})$.
On the other hand, the convolution equivalent tail of $z(1)$ implies a convolution equivalent tail for $u(s)$.
If we further assume that the stationary distribution of $u(s)$ is absolutely continuous, then Theorem~\ref{OUprocess_AI} gives its residual tail dependence function as $\eta(h)=1/(2-e^{-ah})$, which coincides with the limit of its approximation model.
In Figure~\ref{figure::limiteta_OU}, we illustrate this convergence of the residual tail dependence coefficient $\eta$ of the approximating model to the true non-Gaussian OU process for $a=0.2$, $s_1=0$, $T=4$, and three equidistant partitions of the interval $[0,T]$ with mesh length $\Delta=0.4, 0.2, 0.05$. 

The preceding analysis seemingly indicates that the extremal dependence structure is preserved in the limit for convergent (in probability) random vectors.
It is however important to note that in general, convergence in probability or even almost surely does not (!) necessarily imply convergence of the corresponding tail dependence coefficients, as shown in the following counterexample. 

\begin{example}\label{counter_conv}
    Let $X_{1,n} =  X/n + \epsilon_1, X_{2,n} =  X/n + \epsilon_2$, where $X,\epsilon_1,\epsilon_2$ are independent, $X$ has a regularly varying survival function $\Bar{F}_X \in {\text RV}_\rho$, and $\epsilon_i, i=1,2$ are standard normal random variables.
    %We assume $\rho<-2$ to ensure the existence of the second moment of $X$.
    It is straightforward to see that $(X_{1,n}, X_{2,n})$ converges almost surely to the limiting random vector $(\epsilon_1, \epsilon_2)$.
    However, the extremal dependence is clearly not preserved in the limit, since by Proposition 4 in \citet{EOW2019} we know that for any finite $n$, $(X_{1,n}, X_{2,n})$ is asymptotically dependent with tail dependence coefficient $\chi=1$, whilst $(\epsilon_1, \epsilon_2)$ is asymptotically independent with residual tail dependence coefficient $\eta=1/2$.
\end{example}

%Therefore, the underlying reason for the preservation of extremal dependence in the non-Gaussian OU processes needs further investigation and this is beyond the scope of this paper.

%%%%%%%%%%%%%%%%%%%%%%%%%%%%%%%%%%%%%%%% Type G Matern Random Fields %%%%%%%%%%%%%%%%%%%%%%%%%%%%%%%%%%%%%%%%
\subsection{Type G Mat\'ern SPDE Random Fields}\label{TypeG_section}

The second application of our general results is the popular class of type G Mat\'ern SPDE random fields defined as the stationary solution to the SPDE
\begin{equation}\label{TypeGmodel}
    (\kappa^2 - \Delta) ^{\alpha/2} u(\bm{s}) = \dot{\mathcal{M}}(\bm{s}), \quad \bm{s}\in\mathR^d,
\end{equation}
where $\alpha=\nu+d/2$, $d=1,2,\dots$ is the dimension, $\nu>0$ is the smoothness parameter, $\kappa>0$ is the range parameter, $\Delta$ is the Laplacian and $\dot{\mathcal{M}}$ is the so-called type G L\'evy noise \citep{Rosinski1991}.
The solution $u(\bm{s})$ can be expressed as process convolutions
\begin{equation}\label{typeG_integral}
    u(\bm{s}) = \int_{\mathR^d} G(\|\bm{s}-\bm{t}\|) \mathcal{M}({\mathrm d} \bm{t}), \quad \bm{s}\in\mathR^d,
\end{equation}
where $\mathcal{M}$ is the random measure associated with the noise $\dot{\mathcal{M}}$, and $G$ is the Green's function of the differential operator in (\ref{TypeGmodel}) of the form
\begin{equation}\label{typeG_Gfunction}
    G(\|\bm{s}-\bm{t}\|) = \frac{2^{1-\frac{\alpha-d}{2}}}{(4\pi)^{d/2}\Gamma(\alpha/2)\kappa^{\alpha-d}} (\kappa\|\bm{s}-\bm{t}\|)^{\frac{\alpha-d}{2}} K_{\frac{\alpha-d}{2}}(\kappa\|\bm{s}-\bm{t}\|),
\end{equation}
with $\Gamma$ being the gamma function and $K$ being the modified Bessel function of the second kind.

Notably, for the important class of type G Mat\'ern SPDE random fields, the function $G$ is absolutely continuous and monotonically decreasing (see Proposition \ref{Aux_result_3} in Appendix~\ref{Appendix_aux_results}).
This implies that Assumption~\ref{assumptionB1} is satisfied.
Moreover, Assumption~\ref{assumptionB2} is satisfied for the NIG Mat\'ern SPDE fields and variance Gamma Mat\'ern SPDE fields considered in \citet{Bolin2014}, \citet{WallinBolin2015}, and \citet{BolinWallin2020}.
Hence, if we consider its integral approximation of the form~\eqref{mov_ave_approx}, Theorem~\ref{TypeG_limit_eta_prop} gives the limiting residual tail dependence function for convex $G$.

We remark that for $d=2$, the specific Green's function $G$ in \eqref{typeG_Gfunction} is convex only when $\alpha \leq 3$ and $G(0)$ is bounded when $\alpha > 2$.
This implies that when $\alpha=2$, which is the case commonly used in practice, if the constructed mesh is very fine, then the resulting discrete approximation model would have $\eta \approx 1/2$ (near independence) between all pairs of locations, regardless of the distance between them.
From a practical point of view, this indicates that the case $\alpha=2$ might not be suitable for modeling extremal dependence, whilst $\alpha=3$ might be more useful.
We also remark that the frequently used finite element approximation is in the form of linear transformations as well, and more details are given in Appendix~\ref{Appendix_FEM}.

\begin{figure}[!t]
    \centering
    \includegraphics[width=0.8\textwidth]{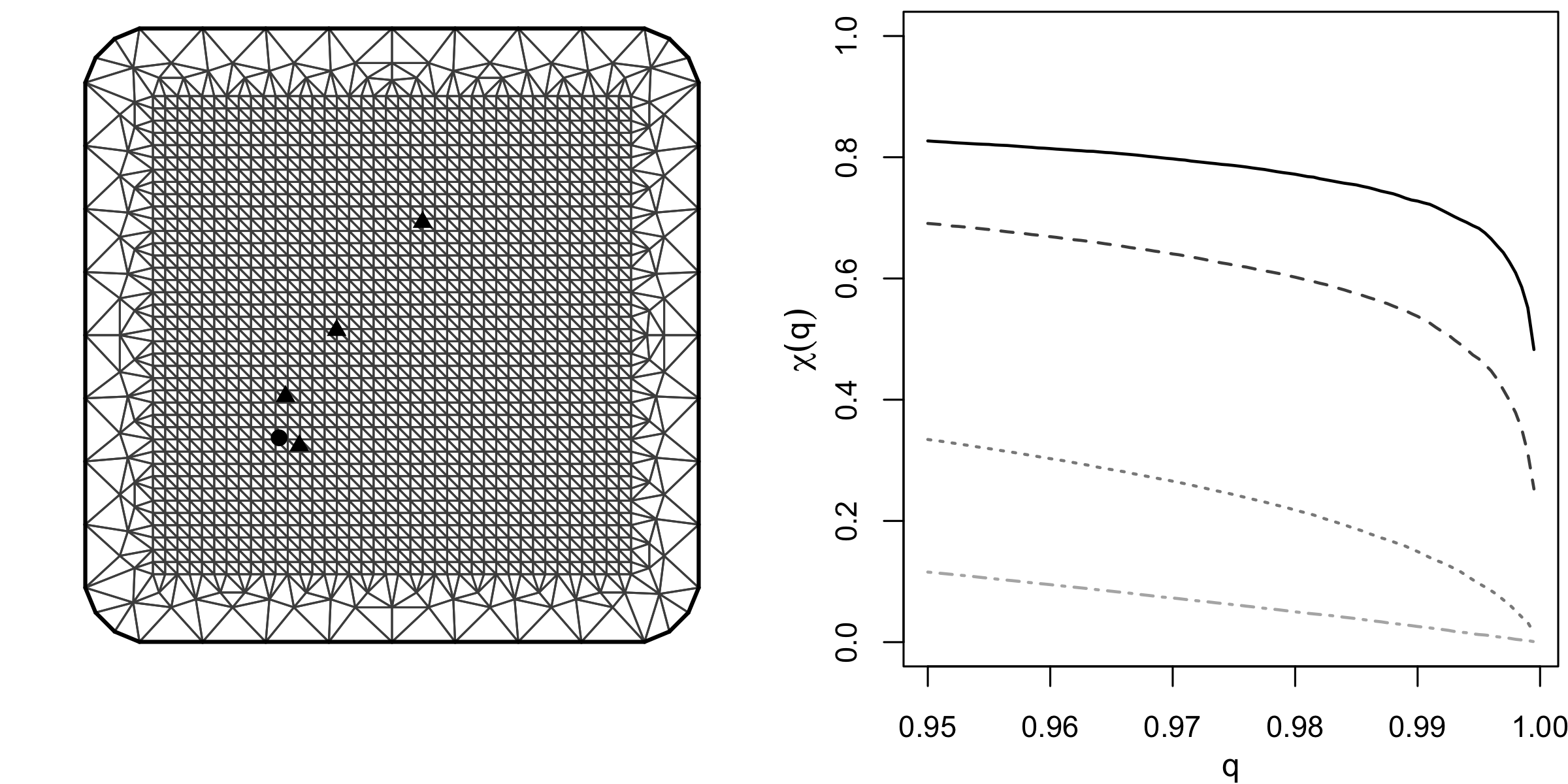}
    \caption{Left panel: constructed mesh and selected pairs of sites with different distances (triangle points to the round point) in the first simulation study in Section~\ref{TypeG_section}; right panel: empirical tail dependence coefficient $\chi(q)$ for different probability levels $q$, of the finite element approximation of the NIG Mat\'ern SPDE model. Four pairs of locations with different distances are plotted, with the lines ordered from top to bottom for increasing distances.}
    \label{figure::NIG_chi_u}
\end{figure}

We conduct simulations to illustrate our results.
We consider the NIG Mat\'ern SPDE model with range parameter $\kappa=2$, smoothness parameter $\alpha=2$, and NIG noise location parameter $\mu=-1$, skewness parameter $\gamma=1$, and shape parameters $\psi=\tau=1$. 
We first consider its finite element approximation and examine the convergence of the pre-asymptotic tail dependence coefficient $\chi(q)$ as $q\rightarrow 1$.
We randomly select 100 sites in the unit square and consider a fine mesh constructed based on a $1600$-node lattice with outer extensions; see the left panel of Figure~\ref{figure::NIG_chi_u}.
Then we simulate $5*10^7$ observations at each site, and compute the empirical tail dependence coefficient $\chi(q)$ with respect to different probability levels $u$.
Figure~\ref{figure::NIG_chi_u} (right panel) shows $\chi(q)$ for four different pairs of locations with different distances.
The plot indicates that the two pairs with longer distances are asymptotically independent and also depicts the decay rate of $\chi(q)$ as $q\rightarrow 1$.
The other two pairs with shorter distances are also asymptotically independent by Proposition~\ref{AInY}, but their corresponding $\chi(q)$ decays at a much slower rate and much more simulations are needed to show that its limit is zero.

We now compare the integral approximation with the finite element approximation and examine the effect of the smoothness parameter $\alpha$ on the extremal dependence structure of the approximation models.
We choose the same SPDE range parameter, NIG noise parameters and the fine mesh constructed based on $1600$ lattice nodes as in the first simulation, and consider more sites, namely 225 randomly selected sites in the unit square, and smoothness parameter $\alpha=2,\dots,5$.
We numerically compute the residual tail dependence coefficient $\eta$ of all pairs of sites for both approximations using the formula from Proposition \ref{AInY}.

\begin{figure}[!t]
    \centering
    \includegraphics[width=\textwidth]{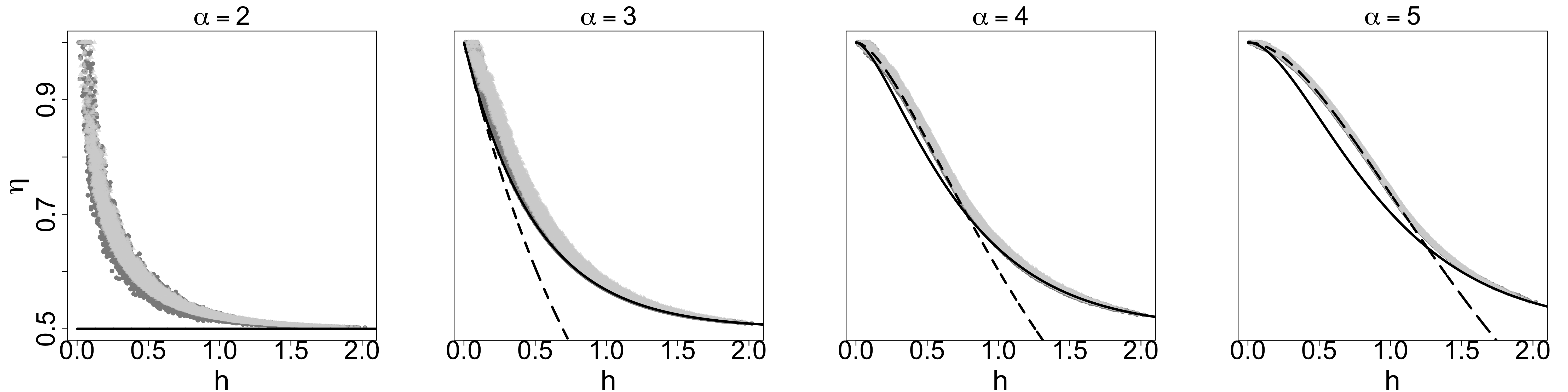}
    \caption{Coefficients $\eta$ for the finite element approximation (light-grey triangle points) and integral approximation (dark-grey round points) of the type G Ma\'ern SPDE model based on a fine mesh, and the functions $1/2+G(h)/\{2G(0)\}$ (black solid line) and $G(h/2)/G(0)$ (black dashed line).}
    \label{figure::NIG_alpha_effect}
\end{figure}

Figure \ref{figure::NIG_alpha_effect} depicts $\eta$ against the distance between all pairs of sites.
A first observation is that with such a fine mesh the difference between the integral approximation and the finite element approximation is negligible.
When $\alpha=2$, i.e., the function $G$ in~\eqref{typeG_Gfunction} is convex but unbounded at $0$, the residual tail dependence function $\eta(h)$ is still fairly far from its limiting function (as the mesh size tends to $0$), which is $\eta(h) \equiv 1/2$.
%seems to have converged to its limiting function, but this is not the case.
%The proof of Theorem~\ref{TypeG_limit_eta_prop} implies that the function depicted by the points in Figure~\ref{figure::NIG_alpha_effect} depends on the mesh size (the shortest distance between the mesh nodes and the observational sites), and this function converges to $1/2$ as the mesh size decreases to $0$.
On the other hand, when $\alpha=3$, i.e., $G$ is convex and bounded at $0$, then the function $\eta(h)$ of both approximation models has almost converged to its limiting function.
While the cases $\alpha=4,5$ are not covered by our theory since $G$ is nonconvex, the $\eta$ values of the approximation models seem to converge to our conjectured limiting function, namely $\eta(h) = {\max\{1/2+G(h)/\{2G(0)\}, G(h/2)/G(0)\}}$. This provides numerical evidence for our conjecture.

\section{Conclusion}
Linear transformations of a random vector with independent components are classical constructions to capture complex correlation dependence structures due to their simplicity and analytical tractability. 
For instance, the approximation models used in practice in the well-known SPDE approach \citep{Lindgren2011,Bolin2014} are of this form.
In this paper, we derived the first results on the extremal dependence structure induced by such constructions when the independent components have exponential tails.
These general results are leveraged to study the extremal dependence structure of moving average processes driven by exponential-tailed L\'evy noise.
In particular, the classical exponential-tailed non-Gaussian OU processes are shown to be asymptotically independent, but with a different residual tail dependence function than their Gaussian counterpart.
As for the type G Mat\'ern random fields, or more general moving average processes, we have shown that under certain assumptions on the kernel function and the noise process, the integral approximation is asymptotically independent when the mesh is fine enough and the limiting residual tail dependence function is derived.

In terms of statistical modeling, linear transformations of exponential-tailed random vectors have the potential to bridge asymptotic dependence and independence in a tractable way, and models with such features are in pursuit in the extremes community \citep{NoldeZhou2021}. 
In fact, this desirable property distinguishes them from other marginal distributions. For instance, linear combinations of heavy-tailed random variables can only result in asymptotic dependence or complete independence. On the other hand,  linear transformations of Gaussian random vectors exclusively yield asymptotic independence or complete dependence. Therefore an interesting and natural question is whether exponential tails are the only ones that can exhibit both extremal dependence classes under linear transformations. 

\section*{Acknowledgement}
We are grateful for Jonathan Tawn's comments on an earlier version of the manuscript.
Zhongwei Zhang and Rapha\"el Huser were partly supported by funding from the King Abdullah University of Science and Technology (KAUST) Office of Sponsored Research (OSR) under Award No. OSR-CRG2020-4394. 

%\section*{Supplementary material}
%\label{SM}
%All proofs, additional simulation studies, auxiliary theoretical results, and more details about the finite element approximation of the type G Mat\'ern SPDE fields are given in the Supplementary material.

\vspace*{-10pt}

%------------------------------------------------------------------------------------
%	                 Appendix
%------------------------------------------------------------------------------------
%\newpage
%\appendix

\begin{appendices}
%%%%%%%%%%%%%%%%%%%%%%%%%%%%%%%%%%%%% Proofs
\vspace*{-10pt}
\section{Lemmas and Proofs for Section \ref{dissection}} \label{Appendix_proofsfor3}

Before proving Proposition \ref{ADnY}, we first present three useful lemmas as follows.
The first one concerns the asymptotic expansion of the quantile of an exponential-tailed distribution function, the second states how much the quantile changes in the limit when an exponential-tailed random variable convolves an independent lighter-tailed random variable, and the third one can be thought of as a more generalized definition of exponential-tailed functions.

\begin{lemma}\label{quantexpotail}
Let $Y$ be a random variable with distribution function $F_Y$ such that $\Bar{F}_Y \in \mathscr{L}_{\beta}, \beta>0$, then we have the asymptotic expansion of its quantile function
\[
F^{-1}_Y(u) \sim -\frac{\log(1-u)}{\beta}, \quad u\uparrow 1.
\]
\end{lemma}
\begin{proof}
%\noindent\textit{Proof of Lemma \ref{quantexpotail}}.
Using the representation (\ref{expotailrepre}), we have
\begin{equation*}
{\bar F}_Y(x) = a(x)\exp\Big\{-\int_0^x \beta(v) {\mathrm d}v  \Big\},
\end{equation*}
where $a(x) \rightarrow a \in (0,\infty)$ and $\beta(x) \rightarrow \beta$ as $x\rightarrow\infty$.
Now consider the function ${y(u)=-\log(1-u)/\beta}$, ${0<u<1}$.
We have
\[
y(F_Y(x)) = -\frac{\log\{\bar{F}_Y(x)\}}{\beta} = \frac{-\log\{a(x)\}+\int_0^x \beta(v) {\mathrm d}v}{\beta} \sim x, \quad x\rightarrow\infty.
\]
Furthermore,
\[
y(F_Y(F_Y^{-1}(u))) = -\frac{\log\{\Bar{F}_Y(F_Y^{-1}(u)\}}{\beta} \sim -\frac{\log(1-u)}{\beta}, \quad u\uparrow 1,
\]
where the last step is due to the fact that $\Bar{F}(F^{-1}(u))\sim 1-u$ and that the logarithm function is slowly varying and thus preserves asymptotic equivalence (see Proposition 2.6 (iii) of \citet{Resnick2007} or Section 3.4 of \citet{Buldygin2018} for a more detailed treatment).
Therefore, combining these two above results gives $F^{-1}_Y(u) \sim -\frac{\log(1-u)}{\beta}$ as $u\uparrow 1$, and the proof is complete. %\hfill \qedsymbol
\end{proof}

%%%%%%%%
\begin{lemma}\label{quantilediff}
Let $Y$ be a random variable with distribution function $F_Y$ such that $\Bar{F}_Y \in\mathscr{L}_\beta, \beta>0$.
Let $Z_1, Z_2$ be random variables independent of $Y$ and have moment generating functions $M_{Z_i}$ such that $M_{Z_i}(\gamma)<\infty$ for some $\gamma>\beta$.
Then for $X_1 = Y + Z_1$, $X_2 = Y + Z_2$, we have
\begin{align*}
    c_1 &:= \lim_{u\uparrow 1} \{F_{X_1}^{-1}(u) - F_{Y}^{-1}(u)\} = \frac{1}{\beta}\log\{M_{Z_1}(\beta)\}, \\
    c_2 &:= \lim_{u\uparrow 1} \{F_{X_2}^{-1}(u) - F_{Y}^{-1}(u)\} = \frac{1}{\beta}\log\{M_{Z_2}(\beta)\}, \\
    c_3 &:= \lim_{u\uparrow 1} \{F_{X_2}^{-1}(u) - F_{X_1}^{-1}(u)\} = c_2 - c_1,
\end{align*}
where $F_{X_1}, F_{X_2}$ are the distribution functions of $X_1$ and $X_2$, respectively.
\end{lemma}
%\noindent\textit{Proof of Lemma \ref{quantilediff}}.
\begin{proof}
By Lemma 1 in \citet{Cline1986} or Proposition 3 in \citet{Breiman1965}, we have that $\Bar{F}_{X_1}, \Bar{F}_{X_2}\in \mathscr{L}_{\beta}$ and
\[
\lim_{x\rightarrow\infty}\frac{{\Bar F}_{X_1}(x)}{{\Bar F}_Y (x)} = M_{Z_1}(\beta), \quad \lim_{x\rightarrow\infty}\frac{{\Bar F}_{X_2}(x)}{{\Bar F}_Y (x)} = M_{Z_1}(\beta).
\]
Using the representation (\ref{expotailrepre}), we have
\[
{\bar F}_{X_1}(x) \sim M_{Z_1}(\beta){\Bar F}_Y (x) \sim M_{Z_1}(\beta) a(x)\exp\Big\{-\int_0^x \beta(v) {\mathrm d}v  \Big\}, \quad x\rightarrow\infty.
\]
Since $\bar{F}_Y(F_Y^{-1}(u)) \sim 1-u \sim \bar{F}_{X_1}(F_{X_1}^{-1}(u))$ as $u$ tends to $1$, we have that as $u\uparrow 1$,
\[
a(F_Y^{-1}(u)) \exp\Big\{-\int_0^{F_Y^{-1}(u)} \beta(v) {\mathrm d}v  \Big\} \sim M_{Z_1}(\beta) a(F_{X_1}^{-1}(u)) \exp\Big\{-\int_0^{F_{X_1}^{-1}(u)} \beta(v) {\mathrm d}v  \Big\}.
\]
As $\Bar{F}_Y, \Bar{F}_{X_1}\in\mathscr{L}_{\beta}$, using Lemma \ref{quantexpotail} we know that 
\[
F_Y^{-1}(u) \sim F_{X_1}^{-1}(u) \sim -\frac{\log(1-u)}{\beta}, \quad u\uparrow 1.
\]
The fact $a(x) \rightarrow a\in(0,\infty)$ implies that the function $a(\cdot)$ is slowly varying and thus it preserves asymptotic equivalence.
That is, $a(F_Y^{-1}(u)) \sim a(F_{X_1}^{-1}(u))$, as $u\uparrow 1$.
Hence, 
\[
M_{Z_1}(\beta) \exp\Big\{-\int_{F_Y^{-1}(u)}^{F_{X_1}^{-1}(u)} \beta(v) {\mathrm d}v  \Big\} \rightarrow 1, \quad u\uparrow 1.
\]
By dominated convergence, we have
\[
M_{Z_1}(\beta)  \exp[-\beta \lim_{u\uparrow 1}\{F_{X_1}^{-1}(u) - F_Y^{-1}(u)\} ] = 1.
\]
Therefore, rearranging the terms gives
\[
c_1 = \lim_{u\uparrow 1} \{ F_{X_1}^{-1}(u) - F_{Y}^{-1}(u) \} =  \frac{1}{\beta}\log\{M_{Z_1}(\beta)\}.
\]
Similarly, we get $c_2 = \log\{M_{Z_2}(\beta)\}/\beta$, and 
\[
c_3 = \lim_{u\uparrow 1} \{F_{X_2}^{-1}(u) - F_{X_1}^{-1}(u)\} = \lim_{u\uparrow 1} [\{F_{X_2}^{-1}(u) - F_{Y}^{-1}(u)\} - \{F_{X_1}^{-1}(u) - F_{Y}^{-1}(u)\} ] = c_2 - c_1.
\] 
%\hfill \qedsymbol
\end{proof}

%%%%%%%%
\begin{lemma}\label{taildiff}
Let $F$ be a distribution function such that $\Bar{F} \in \mathscr{L}_{\beta}, \beta\geq 0$.
Suppose $g_1, g_2$ are two real-valued functions on $\mathR$ that satisfy $g_1(u) \rightarrow \infty$, $g_2(u) \rightarrow \infty$, $g_1(u)/g_2(u) \rightarrow 1$, and ${g_1(u)-g_2(u) \rightarrow g}$ with ${g \in (-\infty, \infty)}$ as $u\rightarrow u_0$, then
\[
\lim_{u\rightarrow u_0}\frac{{\bar F}(g_1(u))}{{\bar F}(g_2(u))} = \exp(-g\beta).
\]
\end{lemma}
\begin{proof}
%\noindent\textit{Proof of Lemma \ref{taildiff}}.
Since $\Bar{F} \in \mathscr{L}_{\beta}$, we know that ${\bar F}(t)$ has the representation (\ref{expotailrepre}).
Using the same argument as in Lemma \ref{quantilediff},  we know that $a(\cdot)$ preserves asymptotic equivalence.
Thus, $a(g_1(u)) \sim a(g_2(u)), u\rightarrow u_0$ and we further have 
\begin{align*}
    \lim_{u\rightarrow u_0}\frac{{\bar F}(g_1(u))}{{\bar F}(g_2(u))} &= \lim_{u\rightarrow u_0} \frac{a(g_1(u))}{a(g_2(u))} \exp\Big\{-\int_{g_2(u)}^{g_1(u)} \beta(v) {\mathrm d}v  \Big\} \\
    &= \lim_{u\rightarrow u_0} \exp\Big\{-\int_{g_2(u)}^{g_1(u)} \beta(v) {\mathrm d}v  \Big\} \\
    &= \exp(-g\beta),
\end{align*}
where the last equality holds by dominated convergence theorem.
This completes the proof. 
%\hfill \qedsymbol
\end{proof}

%%%%%%%%
Now we are ready to prove Proposition \ref{ADnY}. 
\begin{proof}[Proof of Proposition \ref{ADnY}]
Note that extremal dependence is a copula property, i.e., it is invariant to strictly increasing marginal transformations.
This implies that, ${\Tilde{\bm{X}}=(\Tilde{X}_1, \Tilde{X}_2)^\top}$ with ${\Tilde{X}_1=X_1/a_{1 \text{max}}}$, ${\Tilde{X}_2=X_2/a_{2 \text{max}}}$, will have the same extremal dependence structure as $\bm{X}$.
Hence, without losing generality we assume that $a_{1 \text{max}}=a_{2 \text{max}}=1$ and $0\leq a_{ji} \leq 1, j=1,2, i=1,\dots,n$.
We thus have ${Z_1 = \sum_{i\notin I_{\text{max}}} a_{1i}Y_i}$, ${Z_2 = \sum_{i\notin I_{\text{max}}} a_{2i}Y_i}$.
Denote ${Z_c = \sum_{i\in I_{\text{max}}} Y_i}$, then ${X_1=Z_c+Z_1}$ and ${X_2=Z_c+Z_2}$.

In order to derive the tail dependence coefficient $\chi$,  we need to compare the following probability with $1-u$ as $u\uparrow 1$ 
\begin{align*}
    I &:= \Pr\big(X_1 > F_{X_1}^{-1}(u), X_2 > F_{X_2}^{-1}(u)\big) \\
    &= \Pr\big(Z_c + Z_1 > F_{X_1}^{-1}(u), Z_c + Z_2 > F_{X_2}^{-1}(u)\big) \\
    &= \Pr\big(Z_c + Z_1 > F_{X_1}^{-1}(u), Z_c + Z_2 > F_{X_2}^{-1}(u), Z_2 -Z_1 \leq F_{X_2}^{-1}(u) - F_{X_1}^{-1}(u) \big) + \\
    &\quad\quad \Pr\big(Z_c + Z_1 > F_{X_1}^{-1}(u), Z_c + Z_2 > F_{X_2}^{-1}(u), Z_2 -Z_1 > F_{X_2}^{-1}(u) - F_{X_1}^{-1}(u) \big).
\end{align*}
For $Z_j = \sum_{i\notin I_{\text{max}}} a_{ji}Y_i, j=1,2$, if $a_{ji}=0$ for all $i\notin I_{\text{max}}$, then $Z_j=0$ and clearly its moment generating function satisfies $M_{Z_j}(t)=1<\infty$.
If there is at least one $a_{ji} >0$ for some $i\notin I_{\text{max}}$, then by Theorem 3 in \citet{Embrechts1980}, we know that $Z_c$ has exponential-tailed distribution function with index $\beta$, and $Z_j$ has the distribution function $F_{Z_j}$ with $\Bar{F}_{Z_j}\in\mathscr{L}_{\beta/\max_{i\notin I_{\text{max}}} a_{ji}}$.
Hence, $Z_j$ has a lighter exponential tail than $Z_c$ and there exists $\gamma>\beta$ such that $M_{Z_i}(\gamma)<\infty$.
Therefore, the conditions in Lemma \ref{quantilediff} are satisfied.

Since $Z_c + Z_2 > F_{X_2}^{-1}(u)$ and $Z_2 -Z_1 \leq F_{X_2}^{-1}(u) - F_{X_1}^{-1}(u)$ imply $Z_c + Z_1 > F_{X_1}^{-1}(u)$, we have that
\[
\{ Z_c + Z_2 > F_{X_2}^{-1}(u), Z_2 -Z_1 \leq F_{X_2}^{-1}(u) - F_{X_1}^{-1}(u)\} \subseteq \{Z_c + Z_1 > F_{X_1}^{-1}(u) \}.
\]
Hence, 
\begin{align*}
     &\lim_{u\uparrow 1} \frac{\Pr\big(Z_c + Z_1 > F_{X_1}^{-1}(u), Z_c + Z_2 > F_{X_2}^{-1}(u), Z_2 -Z_1 \leq F_{X_2}^{-1}(u) - F_{X_1}^{-1}(u) \big)}{\Pr\big(Z_c>F_{Z_c}^{-1}(u)\big)} \\
    =& \lim_{u\uparrow 1} \frac{\Pr\big(Z_c + Z_2 > F_{X_2}^{-1}(u), Z_2 -Z_1 \leq F_{X_2}^{-1}(u) - F_{X_1}^{-1}(u) \big)}{\Pr\big(Z_c>F_{Z_c}^{-1}(u)\big)} \\
    =& \lim_{u\uparrow 1} \int\int_{z_2 - z_1 \leq F_{X_2}^{-1}(u) - F_{X_1}^{-1}(u)} \frac{\Pr\big(Z_c > F_{X_2}^{-1}(u)-z_2\big)}{\Pr\big(Z_c>F_{Z_c}^{-1}(u)\big)} F_{\bm{Z}}({\mathrm d}z_1,{\mathrm d}z_2) \\
    =& \int\int_{z_2 - z_1 \leq c_3} \exp(\beta z_2 - c_2)  F_{\bm{Z}}({\mathrm d}z_1,{\mathrm d}z_2),
\end{align*}
where the last equality holds due to Lemma \ref{quantilediff}, Lemma \ref{taildiff}, and dominated convergence theorem.
Similarly, we have 
\begin{align*}
     &\lim_{u\uparrow 1} \frac{\Pr\big(Z_c + Z_1 > F_{X_1}^{-1}(u), Z_c + Z_2 > F_{X_2}^{-1}(u), Z_2 -Z_1 > F_{X_2}^{-1}(u) - F_{X_1}^{-1}(u) \big)}{\Pr\big(Z_c>F_{Z_c}^{-1}(u)\big)} \\
    =& \int\int_{z_2 - z_1 > c_3} \exp(\beta z_1 - c_1)  F_{\bm{Z}}({\mathrm d}z_1,{\mathrm d}z_2).
\end{align*}

Therefore, 
\begin{align*}
    \chi &= \lim_{u\uparrow 1} \frac{I}{1-u} \\
    &= \lim_{u\uparrow 1} \frac{\Pr\big(Z_c + Z_1 > F_{X_1}^{-1}(u), Z_c + Z_2 > F_{X_2}^{-1}(u), Z_2 -Z_1 \leq F_{X_2}^{-1}(u) - F_{X_1}^{-1}(u) \big)}{\Pr\big(Z_c>F_{Z_c}^{-1}(u)\big)} + \\
    & \quad\quad \lim_{u\uparrow 1} \frac{\Pr\big(Z_c + Z_1 > F_{X_1}^{-1}(u), Z_c + Z_2 > F_{X_2}^{-1}(u), Z_2 -Z_1 > F_{X_2}^{-1}(u) - F_{X_1}^{-1}(u) \big)}{\Pr\big(Z_c>F_{Z_c}^{-1}(u)\big)} \\
    &= \int\int_{z_2 - z_1 \leq c_3} \exp(\beta z_2 - c_2)  F_{\bm{Z}}({\mathrm d},z_1 {\mathrm d}z_2) + \int\int_{z_2 - z_1 > c_3} \exp(\beta z_1 - c_1)  F_{\bm{Z}}({\mathrm d}z_1,{\mathrm d}z_2) \\
    &= \frac{\int\int_{z_2 - z_1 > c_3} \exp(\beta z_1)  F_{\bm{Z}}({\mathrm d}z_1,{\mathrm d}z_2)}{M_{Z_1}(\beta)} +
    \frac{\int\int_{z_2 - z_1 \leq c_3} \exp(\beta z_2)  F_{\bm{Z}}({\mathrm d}z_1,{\mathrm d}z_2)}{M_{Z_2}(\beta)} \\
    &= \E[\min\{\exp(\beta Z_1)/M_{Z_1}(\beta), \exp(\beta Z_2)/M_{Z_2}(\beta) \}].
\end{align*} 
\end{proof}

%%%%%%%%
We now prove Proposition \ref{ADtwoY_limit}.
\begin{proof}[Proof of Proposition \ref{ADtwoY_limit}]
For $Y_1 \sim {\text{GH}}(\lambda, \tau, \psi, \mu=0, \gamma=0)$, from Section \ref{ExpoIntro} we know that its density function $f_{Y_1}(y)$ has asymptotic expansion
\begin{equation*}
    f_{Y_1}(y) \sim c_0 y^{\lambda-1}e^{-\sqrt{\psi} y}, \quad y\rightarrow \infty,
\end{equation*}
where $c_0>0$ is a constant.
It can be easily seen that its moment generating function $M_{Y_1}(t)$ evaluated at point $t=\sqrt{\psi}$ is finite if $\lambda<0$, but is infinite if $\lambda\geq0$.
Hence, if $\lambda<0$, we have 
$$
\lim_{a_{22}\uparrow 1} \chi =  \frac{\int_{-\infty}^{c^*} \exp(\sqrt{\psi}y)F_{Y_1}({\mathrm d}y)}{M_{Y_1}(\sqrt{\psi})} + 
\frac{\int^{\infty}_{c^*} \exp(a_{12}\sqrt{\psi}y)  F_{Y_1}({\mathrm d}y)}{M_{Y_1}(a_{12}\sqrt{\psi})}.
$$
    
If $\lambda\geq0$, we know that $\lim_{a_{22}\uparrow 1} M_{Y_1}(a_{22}\sqrt{\psi}) = \infty$ and thus $\lim_{a_{22}\uparrow 1} c = \infty$.
Hence, the second term in the expression of $\chi$ in Example \ref{ADtwoY} vanishes.
We now show the first term tends to zero as $a_{22}\uparrow 1$.
Since $f_{Y_1}(y) \sim c_0 y^{\lambda-1}e^{-\sqrt{\psi} y}, y\rightarrow \infty$, we have for any $\varepsilon >0$, there exists $N >0$ such that ${f_{Y_1}(y) < (1+\varepsilon)c_0 y^{\lambda-1}e^{-\sqrt{\psi} y}}$ for every $y>N$.
Hence, for $\lambda>0$
\begin{align}
    &\quad\frac{\int_{-\infty}^{c/(a_{22}-a_{12})} \exp(a_{22}\sqrt{\psi}y)F_{Y_1}({\mathrm d}y)}{M_{Y_1}(a_{22}\sqrt{\psi})}  \nonumber\\
    &= \frac{\int_{-\infty}^{c/(a_{22}-a_{12})} \exp(a_{22}\sqrt{\psi}y)f_{Y_1}(y){\mathrm d}y}{M_{Y_1}(a_{22}\sqrt{\psi})} \nonumber\\
    &< \frac{\int_{-\infty}^{N} \exp(a_{22}\sqrt{\psi}y)f_{Y_1}(y){\mathrm d}y}{M_{Y_1}(a_{22}\sqrt{\psi})} + \frac{\int_{N}^{c/(a_{22}-a_{12})} (1+\varepsilon)cy^{\lambda-1}\exp\{-(1-a_{22})\sqrt{\psi}y\}{\mathrm d}y}{M_{Y_1}(a_{22}\sqrt{\psi})} \nonumber\\
    &< \frac{\exp(a_{22}\sqrt{\psi}N) \int_{-\infty}^{N} f_{Y_1}(y){\mathrm d}y}{M_{Y_1}(a_{22}\sqrt{\psi})} + \frac{\int_{N}^{c/(a_{22}-a_{12})} (1+\varepsilon)cy^{\lambda-1}{\mathrm d}y}{M_{Y_1}(a_{22}\sqrt{\psi})} \nonumber\\
    &= \frac{\exp(a_{22}\sqrt{\psi}N) F_{Y_1}(N)}{M_{Y_1}(a_{22}\sqrt{\psi})} + \frac{ (1+\varepsilon)c [\{c/(a_{22}-a_{12})\}^\lambda - N^\lambda]/\lambda}{M_{Y_1}(a_{22}\sqrt{\psi})} \label{lambda>0}.
\end{align}
Substituting $c$ by $c = \frac{\log\{M_{Y_1}(a_{22}\sqrt{\psi})\} - \log\{M_{Y_1}(a_{12}\sqrt{\psi})\}}{\sqrt{\psi}}$, one can easily see that (\ref{lambda>0}) tends to zero as $a_{22}\uparrow 1$.
Hence, we have 
\begin{equation*}
    0 \leq \liminf_{a_{22}\uparrow 1} \frac{\int_{-\infty}^{c/(a_{22}-a_{12})} \exp(a_{22}\sqrt{\psi}y)F_{Y_1}({\mathrm d}y)}{M_{Y_1}(a_{22}\sqrt{\psi})} \leq \limsup_{a_{22}\uparrow 1} \frac{\int_{-\infty}^{c/(a_{22}-a_{12})} \exp(a_{22}\sqrt{\psi}y)F_{Y_1}({\mathrm d}y)}{M_{Y_1}(a_{22}\sqrt{\psi})} = 0,
\end{equation*}
and thus $\lim_{a_{22}\uparrow 1} \frac{\int_{-\infty}^{c/(a_{22}-a_{12})} \exp(a_{22}\sqrt{\psi}y)F_{Y_1}({\mathrm d}y)}{M_{Y_1}(a_{22}\sqrt{\psi})} = 0$.
Note that for $\lambda=0$, we only need to change the term $(c^\lambda - N^\lambda)/\lambda$ in (\ref{lambda>0}) to $(\log c - \log N)$ and the same conclusion holds.
Therefore, for $\lambda\geq 0$, we have $\lim_{a_{22}\uparrow 1} \chi = 0$.
This completes the proof.
\end{proof}

%%%%%%%%%%%%%%%%%% AI case 
%%%%%%%%
We now present four useful lemmas which are needed for proving Proposition \ref{AInY}.
The former three concern the analytical solution to a certain convex optimization problem, while the last one presents some properties of the density function of an exponential-tailed distribution function.

\begin{lemma}\label{optimization}
For $n\geq 3$, let $g_n(x_1, x_2)$ be the minimizer of the following unconstrained optimization problem
\begin{align*}
g_n(x_1, x_2) &= \min_{x_3,\dots,x_n} \Big|\frac{a_{22}x_1-a_{12}x_2+(a_{12}a_{23}-a_{13}a_{22})x_3+\cdots+(a_{12}a_{2n}-a_{1n}a_{22})x_n}{a_{11}a_{22}-a_{12}a_{21}} \Big| + \\
&\Big|\frac{-a_{21}x_1+a_{11}x_2+(a_{21}a_{13}-a_{23}a_{11})x_3+\cdots+(a_{21}a_{1n}-a_{2n}a_{11})x_n}{a_{11}a_{22}-a_{12}a_{21}} \Big| + \sum_{i=3}^n |x_i|
\end{align*}
with $a_{11}a_{22}-a_{12}a_{21}\neq 0$, and $a_{1i}^2+a_{2i}^2\neq 0, i=1,\dots,n$. 
Then $g_n(x_1, x_2)$ can be expressed explicitly as
\[
g_n(x_1, x_2) = \min_{\substack{i,j=1,\dots,n \\ i\neq j}} \frac{|a_{2i}x_1-a_{1i}x_2|+|a_{2j}x_1-a_{1j}x_2|}{|a_{2i}a_{1j}-a_{1i}a_{2j}|}.
\]
\end{lemma}
%\noindent\textit{Proof of Lemma \ref{optimization}}.
\begin{proof}
We prove the statement by induction.
Denote the objective function in $g_n(x_1, x_2)$ as $h_n(x_1,\dots,x_n)$.
For $n=3$, we have
\begin{align*}
g_3(x_1, x_2) &= \min_{x_3} h_3(x_1,x_2,x_3) = \min_{x_3}
\frac{|a_{22}x_1-a_{12}x_2+(a_{12}a_{23}-a_{13}a_{22})x_3|}{|a_{11}a_{22}-a_{12}a_{21}|} + \\
&\quad\quad\quad\quad\quad \frac{|-a_{21}x_1+a_{11}x_2+(a_{21}a_{13}-a_{23}a_{11})x_3|}{|a_{11}a_{22}-a_{12}a_{21}|}+ |x_3|.
\end{align*}
Note that if we fix $x_1, x_2$, then one can partition the real line $\mathR$ into four parts (two rays and two segments) such that the objective function $h_3$ is linear in each part.
Obviously the minimizer cannot be achieved when $x_3=\pm\infty$.
Hence, it must be achieved at the boundary of one of these segments.
That is, if ${a_{12}a_{23}-a_{13}a_{22}\neq 0}$ and ${a_{21}a_{13}-a_{23}a_{11}\neq 0}$, then
\begin{align*}
    g_3(x_1, x_2) = \min \Big\{ h_3(x_1, x_2, 0), h_3\Big(x_1,x_2,\frac{a_{22}x_1-a_{12}x_2}{a_{22}a_{13}-a_{12}a_{23}}\Big), 
    h_3\Big(x_1,x_2,\frac{a_{21}x_1-a_{11}x_2}{a_{21}a_{13}-a_{11}a_{23}}\Big) \Big\}.
\end{align*}
Clearly, $h_3(x_1,x_2,0)=(|a_{22}x_1-a_{12}x_2|+|a_{21}x_1-a_{11}x_2|)/|a_{22}a_{11}-a_{12}a_{21}|$.
Furthermore,
\begin{align*}
    h_3\Big(x_1,x_2,\frac{a_{22}x_1-a_{12}x_2}{a_{22}a_{13}-a_{12}a_{23}}\Big) &= \frac{|-a_{21}x_1+a_{11}x_2+\frac{(a_{22}x_1-a_{12}x_2)(a_{21}a_{13}-a_{23}a_{11})}{a_{22}a_{13}-a_{12}a_{23}}|}{|a_{11}a_{22}-a_{12}a_{21}|} + \Big|\frac{a_{22}x_1-a_{12}x_2}{a_{22}a_{13}-a_{12}a_{23}} \Big| \\
    &= \Big|\frac{a_{23}x_1-a_{13}x_2}{a_{22}a_{13}-a_{12}a_{23}} \Big| + \Big|\frac{a_{22}x_1-a_{12}x_2}{a_{22}a_{13}-a_{12}a_{23}} \Big|.
\end{align*}
Similarly, we have 
\[
h_3\Big(x_1,x_2,\frac{a_{21}x_1-a_{11}x_2}{a_{21}a_{13}-a_{11}a_{23}}\Big) = \Big|\frac{a_{23}x_1-a_{13}x_2}{a_{13}a_{21}-a_{11}a_{23}} \Big| + \Big|\frac{a_{21}x_1-a_{11}x_2}{a_{13}a_{21}-a_{11}a_{23}} \Big|.
\]
Hence, 
\[
g_3(x_1, x_2) = \min_{\substack{i,j=1,2,3 \\ i\neq j}} \frac{|a_{2i}x_1-a_{1i}x_2|+|a_{2j}x_1-a_{1j}x_2|}{|a_{2i}a_{1j}-a_{1i}a_{2j}|}.
\]

If $a_{12}a_{23}-a_{13}a_{22}= 0$ and $a_{21}a_{13}-a_{23}a_{11}= 0$, the above statement clearly also holds.
If only one of them is zero, without loss of generality we assume $a_{12}a_{23}-a_{13}a_{22}\neq 0$ and $a_{21}a_{13}-a_{23}a_{11}= 0$.
Since $a_{11}^2+a_{21}^2\neq 0, a_{13}^2+a_{23}^2\neq 0$, we further assume $a_{11}\neq0, a_{13}\neq 0$ without losing generality.
Then $a_{21}=a_{23}a_{11}/a_{13}$, $g_3(x_1, x_2) = \min \Big\{ h_3(x_1, x_2, 0), h_3\Big(x_1,x_2,\frac{a_{22}x_1-a_{12}x_2}{a_{22}a_{13}-a_{12}a_{23}}\Big) \Big\}$, and
\begin{align*}
    h_3\Big(x_1,x_2,\frac{a_{22}x_1-a_{12}x_2}{a_{22}a_{13}-a_{12}a_{23}}\Big)  &= 
    \frac{|-a_{21}x_1+a_{11}x_2|}{|a_{11}a_{22}-a_{12}a_{21}|} + \Big|\frac{a_{22}x_1-a_{12}x_2}{a_{22}a_{13}-a_{12}a_{23}} \Big| \\
    &= \Big|\frac{a_{23}a_{11}x_1/a_{13}-a_{11}x_2}{a_{11}a_{22}-a_{12}a_{23}a_{11}/a_{13}} \Big| + \Big|\frac{a_{22}x_1-a_{12}x_2}{a_{22}a_{13}-a_{12}a_{23}} \Big| \\
    &= \frac{|a_{23}x_1-a_{13}x_2|+|a_{22}x_1-a_{12}x_2|}{|a_{22}a_{13}-a_{12}a_{23}|}.
\end{align*}
Hence, the statement $g_3(x_1, x_2) = \min_{\substack{i,j=1,2,3 \\ i\neq j}} \frac{|a_{2i}x_1-a_{1i}x_2|+|a_{2j}x_1-a_{1j}x_2|}{|a_{2i}a_{1j}-a_{1i}a_{2j}|}$ also holds in this case.
Above all, we have shown that 
\[
g_3(x_1, x_2) = \min_{\substack{i,j=1,2,3 \\ i\neq j}} \frac{|a_{2i}x_1-a_{1i}x_2|+|a_{2j}x_1-a_{1j}x_2|}{|a_{2i}a_{1j}-a_{1i}a_{2j}|}.
\]

Now we assume that $g_n(x_1, x_2) = \min_{\substack{i,j=1,\dots,n \\ i\neq j}} \frac{|a_{2i}x_1-a_{1i}x_2|+|a_{2j}x_1-a_{1j}x_2|}{|a_{2i}a_{1j}-a_{1i}a_{2j}|}$ holds for some $n\geq 3$ and consider the case $n+1$.
Since there is no constraint in our optimization problem (thus the constraints are independent), we can minimize the function $h_{n+1}$ by first minimizing over some variables, and then minimizing over the remaining variables; see Section 4.1.3 of \citet{Boyd2004} for more details.
Hence,
\begin{align*}
    g_{n+1}(x_1,x_2) = \min_{x_3,\dots,x_{n+1}}h_{n+1}(x_1,\dots,x_{n+1}) = \min_{x_3,\dots,x_n} \min_{x_{n+1}} h_{n+1}(x_1,\dots,x_{n+1}).
\end{align*}
Then using the same arguments as $n=3$, if  ${a_{12}a_{2(n+1)}-a_{1(n+1)}a_{22}\neq 0}$, ${a_{21}a_{1(n+1)}-a_{2(n+1)}a_{11}\neq 0}$, we have
\begin{align*}
    g_{n+1}(x_1,x_2) &= \min_{x_3,\dots,x_n} \min\Big\{ h_{n+1}(x_1, \dots, x_n, 0), \\
    &\quad\quad\quad\quad\quad\quad
    h_{n+1}\Big(x_1,\dots,x_n,\frac{a_{22}x_1-a_{12}x_2+\sum_{i=3}^n (a_{12}a_{2i}-a_{1i}a_{22})x_i}{a_{1(n+1)}a_{22}-a_{12}a_{2(n+1)}}\Big), \\
    &\quad\quad\quad\quad\quad\quad h_{n+1}\Big(x_1,\dots,x_n,\frac{a_{21}x_1-a_{11}x_2+\sum_{i=3}^n (a_{11}a_{2i}-a_{1i}a_{21})x_i}{a_{1(n+1)}a_{21}-a_{11}a_{2(n+1)}}\Big) \Big\} \\
    &= \min\Big\{ \min_{x_3,\dots,x_n} h_{n+1}(x_1, \dots, x_n, 0), \\
    &\quad\quad\quad\quad  \min_{x_3,\dots,x_n}
    h_{n+1}\Big(x_1,\dots,x_n,\frac{a_{22}x_1-a_{12}x_2+\sum_{i=3}^n (a_{12}a_{2i}-a_{1i}a_{22})x_i}{a_{1(n+1)}a_{22}-a_{12}a_{2(n+1)}}\Big), \\
    &\quad\quad\quad\quad  \min_{x_3,\dots,x_n} h_{n+1}\Big(x_1,\dots,x_n,\frac{a_{21}x_1-a_{11}x_2+\sum_{i=3}^n (a_{11}a_{2i}-a_{1i}a_{21})x_i}{a_{1(n+1)}a_{21}-a_{11}a_{2(n+1)}}\Big) \Big\}.
\end{align*}
Note that $\min_{x_3,\dots,x_n} h_{n+1}(x_1, \dots, x_n, 0)=\min_{x_3,\dots,x_n} h_n(x_1, \dots, x_n)=g_n(x_1,x_2)$, and
\begin{align*}
    &\min_{x_3,\dots,x_n}
    h_{n+1}\Big(x_1,\dots,x_n,\frac{a_{22}x_1-a_{12}x_2+\sum_{i=3}^n (a_{12}a_{2i}-a_{1i}a_{22})x_i}{a_{1(n+1)}a_{22}-a_{12}a_{2(n+1)}}\Big) \\
    =& \min_{x_3,\dots,x_n} \Big|\frac{-a_{21}x_1+a_{11}x_2+\sum_{i=1}^n (a_{21}a_{1i}-a_{2i}a_{11})x_i + \frac{a_{21}a_{1(n+1)}-a_{2(n+1)}a_{11}}{a_{1(n+1)}a_{22}-a_{12}a_{2(n+1)}}\{a_{22}x_1-a_{12}x_2+}{a_{11}a_{22}-a_{12}a_{21}} \\
    &\quad\quad \frac{\sum_{i=3}^n (a_{12}a_{2i}-a_{1i}a_{22})x_i\}}{a_{11}a_{22}-a_{12}a_{21}} \Big| + \sum_{i=3}^n |x_i| + \Big| \frac{a_{22}x_1-a_{12}x_2+\sum_{i=3}^n (a_{12}a_{2i}-a_{1i}a_{22})x_i}{a_{1(n+1)}a_{22}-a_{12}a_{2(n+1)}} \Big| \\
    =& \min_{x_3,\dots,x_n}\Big|\frac{a_{2(n+1)}x_1-a_{1(n+1)}x_2+\sum_{i=3}^n (a_{1(n+1)}a_{2i}-a_{1i}a_{2(n+1)})x_i}{a_{1(n+1)}a_{22}-a_{12}a_{2(n+1)}} \Big| + \sum_{i=3}^n |x_i| + \\
    & \quad\quad\quad  \Big| \frac{a_{22}x_1-a_{12}x_2+\sum_{i=3}^n (a_{12}a_{2i}-a_{1i}a_{22})x_i}{a_{1(n+1)}a_{22}-a_{12}a_{2(n+1)}} \Big|,
\end{align*}
which is equal to $g_n(x_1,x_2)$ if we replace the coefficients $a_{11},a_{21}$ in $g_n(x_1,x_2)$ with $a_{1(n+1)}$ and $a_{2(n+1)}$, respectively.
Hence,
\begin{align*}
&\min_{x_3,\dots,x_n}
    h_{n+1}\Big(x_1,\dots,x_n,\frac{a_{22}x_1-a_{12}x_2+\sum_{i=3}^n (a_{12}a_{2i}-a_{1i}a_{22})x_i}{a_{1(n+1)}a_{22}-a_{12}a_{2(n+1)}}\Big) \\
=&\min_{\substack{i,j=2,\dots,n+1 \\ i\neq j}} \frac{|a_{2i}x_1-a_{1i}x_2|+|a_{2j}x_1-a_{1j}x_2|}{|a_{2i}a_{1j}-a_{1i}a_{2j}|}.
\end{align*}
Similarly, we have
\begin{align*}
&\min_{x_3,\dots,x_n} h_{n+1}\Big(x_1,\dots,x_n,\frac{a_{21}x_1-a_{11}x_2+\sum_{i=3}^n (a_{11}a_{2i}-a_{1i}a_{21})x_i}{a_{1(n+1)}a_{21}-a_{11}a_{2(n+1)}}\Big) \\
=& \min_{\substack{i,j=1,3,\dots,n+1 \\ i\neq j}} \frac{|a_{2i}x_1-a_{1i}x_2|+|a_{2j}x_1-a_{1j}x_2|}{|a_{2i}a_{1j}-a_{1i}a_{2j}|}.
\end{align*}
Above all, 
\begin{align*}
    &g_{n+1}(x_1,x_2) = \min\Big\{\min_{\substack{i,j=1,\dots,n \\ i\neq j}} \frac{|a_{2i}x_1-a_{1i}x_2|+|a_{2j}x_1-a_{1j}x_2|}{|a_{2i}a_{1j}-a_{1i}a_{2j}|}, \\ 
    & \min_{\substack{i,j=2,\dots,n+1 \\ i\neq j}} \frac{|a_{2i}x_1-a_{1i}x_2|+|a_{2j}x_1-a_{1j}x_2|}{|a_{2i}a_{1j}-a_{1i}a_{2j}|}, \min_{\substack{i,j=1,3,\dots,n+1 \\ i\neq j}} \frac{|a_{2i}x_1-a_{1i}x_2|+|a_{2j}x_1-a_{1j}x_2|}{|a_{2i}a_{1j}-a_{1i}a_{2j}|} \Big\} \\
    &\quad\quad\quad\quad\quad  =\min_{\substack{i,j=1,\dots,n+1 \\ i\neq j}} \frac{|a_{2i}x_1-a_{1i}x_2|+|a_{2j}x_1-a_{1j}x_2|}{|a_{2i}a_{1j}-a_{1i}a_{2j}|}.
\end{align*}

If  $a_{12}a_{2(n+1)}-a_{1(n+1)}a_{22}=a_{21}a_{1(n+1)}-a_{2(n+1)}a_{11}= 0$ or only one of them is zero, it is straightforward to show that the above statement also holds by using the argument as in the case $n=3$.
Therefore, the proof is complete.  
%\hfill \qedsymbol
\end{proof}

%%%%%%%%
\begin{lemma}\label{optimization12}
Let $a,b,c,d \in [0,\infty)$ with $ad-bc\neq 0$, and $f(x_1,x_2)=\frac{|ax_1-bx_2| + |cx_1-dx_2|}{|ad-bc|}$, then we have 
\[
\min_{x_1\geq 1, x_2\geq 1} f(x_1,x_2) = \min\Big\{\max(1/a,1/b), \max(1/c,1/d), \frac{|a-b| + |c-d|}{|ad-bc|} \Big\}.
\]
\end{lemma}
%\noindent\textit{Proof of Lemma \ref{optimization12}}.
\begin{proof}
Since one can always partition $\mathR^2$ into four subregions such that the function $f(x_1,x_2)$ is linear in each subregion, we know that the minimizer of this function on $[1,\infty)^2$ can only be achieved at the intersection of one subregion and $[1,\infty)^2$.
That is, if $a\leq b, c\leq d$, then
\begin{align*}
    \min_{x_1\geq 1, x_2\geq 1} f(x_1,x_2) &= \min\{f(1,1), f(b/a,1), f(d/c,1) \} \\
    &= \min\Big\{\frac{|a-b| + |c-d|}{|ad-bc|}, 1/a, 1/c \Big\} \\
    &= \min\Big\{\frac{|a-b| + |c-d|}{|ad-bc|}, \max(1/a,1/b), \max(1/c,1/d)\Big\}.
\end{align*}
Similarly, for $b\leq a, c\leq d$, or $a\leq b, d\leq c$, or $b\leq a, d\leq c$, we also obtain
\[
\min_{x_1\geq 1, x_2\geq 1} f(x_1,x_2) = \min\Big\{\frac{|a-b| + |c-d|}{|ad-bc|}, \max(1/a,1/b), \max(1/c,1/d)\Big\}.
\]
Therefore, above all, the conclusion always holds.
%\hfill \qedsymbol
\end{proof}

%%%%%%%%
\begin{lemma}\label{inequalities}
Let $a,b,c,d\in [0,1]$ with $ad-bc\neq 0$, $ab<1$ and $cd<1$, then
\[
  \frac{|a-b|+|c-d|}{|ad-bc|} > 1.
\]
% $\min(1/a,1/b) > \frac{2-a-b}{1-ab} > 1$.
\end{lemma}
%\noindent\textit{Proof of Lemma \ref{inequalities}}.
\begin{proof}
Since $ad-bc\neq 0$, we know that $a-b$ and $c-d$ cannot be both equal to zero.
Hence,
\begin{equation*}
    |ad-bc| = |ad-ac+ac-bc| \leq a|c-d|+c|a-b| < |c-d| + |a-b|,
\end{equation*}
where the last inequality is due to the fact that $ab<1$ and $cd<1$, which means that $a$ and $b$ cannot be both equal to $1$, and $c$ and $d$ cannot be both equal to $1$.
This implies that $\frac{|a-b|+|c-d|}{|ad-bc|} > 1$. 
%\hfill \qedsymbol
\end{proof}

%%%%%%%%
\begin{lemma}\label{expotaildens}
Let $Y$ be a random variable with distribution function $F_Y$ such that $\Bar{F}_Y\in\mathscr{L}_\beta, \beta>0$.
If $F_Y$ is absolutely continuous with density function $f_Y$, then we have the asymptotic relationship
\[
\log \bar{F}_Y(x) \sim \log f_Y(x) \sim -\beta x, \quad x\rightarrow\infty.
\]
\end{lemma}
%\noindent\textit{Proof of Lemma \ref{expotaildens}}.
\begin{proof}
Since $\Bar{F}_Y \in \mathscr{L}_\beta$, using L'H\^opital's rule, we know that $f_Y \in \mathscr{L}_\beta$.
Now suppose that ${g: \mathR \rightarrow \mathR_+}$ is an exponential-tailed function with index $\beta$, i.e., $g\in\mathscr{L}_\beta$.
Using the representation (\ref{expotailrepre}), we have
\begin{equation*}
g(x) = a(x)\exp\Big\{-\int_0^x \beta(v) {\mathrm d}v  \Big\},
\end{equation*}
where $a(x) \rightarrow a \in (0,\infty)$ and $\beta(x) \rightarrow \beta$ as $x\rightarrow\infty$.
The fact $\beta(x) \rightarrow \beta$ implies that $\beta(\cdot)$ is slowly regularly.
By Karamata's theorem (see Theorem 2.1 in \citet{Resnick2007}), we have that 
\[
\int_0^x \beta(v) {\mathrm d}v \sim \beta x, \quad x\rightarrow \infty.
\]
Hence,
\begin{align*}
    \log g(x) = \log a(x) - \int_0^x \beta(v) {\mathrm d}v 
                \sim - \int_0^x \beta(v) {\mathrm d}v
                \sim - \beta x, \quad x\rightarrow \infty.
\end{align*}
Since both $f_Y$ and $\bar{F}_Y$ have exponential tails with index $\beta$, we have the asymptotic expansion ${\log \bar{F}_Y(x) \sim \log f_Y(x) \sim -\beta x}$, ${x\rightarrow\infty}$. 
%\hfill \qedsymbol
\end{proof}

%%%%%%%%
We are now ready to prove Proposition \ref{AInY}. 
\begin{proof}[Proof of Proposition \ref{AInY}]
If~\eqref{AD_case}, then clearly $\eta=1$, which coincides with the result in Proposition~\ref{ADnY}.
Otherwise, without losing generality we assume that $a_{11}=\max_{i\in\{1,\dots,n\}} a_{1i}$ and $a_{22}=\max_{i\in\{1,\dots,n\}} a_{2i}$.
Since extremal dependence is a copula property, $(X_1/a_{11}, X_2/a_{22})$ has the same extremal dependence structure as $\bm{X}$.
Hence, we can further assume that ${a_{11}=a_{22}=1}$ and ${0\leq a_{ji}\leq 1}$, ${j=1,2}$, ${i=1,\dots,n}$.
Since ${\argmax_{i\in\{1,\dots,n\}} a_{1i} \cap \argmax_{i\in\{1,\dots,n\}} a_{2i} =\emptyset}$, we know that ${a_{1i}a_{2i}<1}$.

The strategy for this proof relies on augmenting the model ${X_1 = Y_1 + a_{12}Y_2+\cdots+a_{1n}Y_n}$, ${X_2 = a_{21}Y_1 + Y_2+a_{23}Y_3+\cdots+a_{2n}Y_n}$ in the following way
\begin{empheq}[left=\empheqlbrace]{align*}
    X_1 &= Y_1 + a_{12}Y_2+\cdots+a_{1n}Y_n, \\
    X_2 &= a_{21}Y_1 + Y_2+a_{23}Y_3+\cdots+a_{2n}Y_n, \\
    X_3 &= Y_3, \\
    &\vdots \\
    X_n &= Y_n.
\end{empheq}
Since $Y_i, i=1,\dots,n$ are independent and have common distribution function $\Bar{F}_Y\in\mathscr{L}_\beta$, by Theorem 3 in \citet{Embrechts1980}, we know that $\Bar{F}_{X_i}\in\mathscr{L}_{\beta}, i=1,\dots,n$.
Using Lemma \ref{expotaildens}, we have
\begin{equation*}
   1 - F_{X_i}(x) = e^{-h(x)}, \quad h(x)\sim \beta x, \quad \text{ as } x\rightarrow\infty.
\end{equation*}
Now for the square matrix 
\begin{equation*}
A=\begin{pmatrix} 
1 & a_{12} & a_{13} & \cdots & a_{1n} \\ 
a_{21} & 1 & a_{23} & \cdots & a_{2n} \\
 &  & 1 &  \\
 &  &  & \ddots & \\
 &  &  &   &  1 
\end{pmatrix} \in \mathR^{n\times n},
\end{equation*}
it can be shown that its determinant is $|A| = 1-a_{12}a_{21}$, and its inverse is
\begin{equation*}
A^{-1}= \frac{1}{1-a_{12}a_{21}}
\begin{pmatrix} 
1 & -a_{12} & a_{12}a_{23}-a_{13} & \cdots & a_{12}a_{2n}-a_{1n} \\ 
-a_{21} & 1 & a_{21}a_{13}-a_{23} & \cdots & a_{21}a_{1n}-a_{2n} \\
 &  & 1-a_{12}a_{21} &  \\
 &  &  & \ddots & \\
 &  &  &   &  1-a_{12}a_{21} 
\end{pmatrix}.
\end{equation*}

Denote $\bm{X}=(X_1, X_2)^T$, $\Tilde{\bm{X}}=(X_1, \dots, X_n)^T$ and the probability density function of $\Tilde{\bm{X}}$ as $f_{\Tilde{\bm{X}}}$.
Then we have 
\begin{align*}
    \lim_{t\rightarrow\infty} \frac{-\log f_{\Tilde{\bm{X}}}(t\bm{x})}{h(t)} &= \lim_{t\rightarrow\infty}\frac{\beta t (|\bm{\alpha}_1^T \bm{x}| + \cdots + |\bm{\alpha}_n^T \bm{x}|)} {\beta t} \\
    &= |\bm{\alpha}_1^T \bm{x}| + \cdots + |\bm{\alpha}_n^T \bm{x}| \\
    &=: \Tilde{g}(\bm{x}),
\end{align*}
where $\bm{\alpha}_i, i=1,\dots,n$ are the $i$th row vector of the matrix $A^{-1}$.
By Proposition 2 in \citet{Nolde2022}, a sequence of scaled random samples $N_n=\{\Tilde{\bm{X}}_1/r_n, \dots, \Tilde{\bm{X}}_n/r_n\}$ from $f_{\Tilde{\bm{X}}}$ converges in probability onto a limit set $\Tilde{G}$ with $\Tilde{G}=\{\bm{x}\in\mathR^d: \Tilde{g}(\bm{x}) \leq 1\}$.
Then using Proposition 4 in \citet{Nolde2022} we know that, for $\bm{X}$, which is a two-dimensional subvector of $\Tilde{\bm{X}}$, sample clouds from $\bm{X}$ converge onto the limit set $G=\{x_1,x_2\in\mathR: g(x_1,x_2) \leq 1\}$ with gauge function
\begin{equation*}
    g(x_1,x_2) = \min_{x_3,\dots,x_n} \Tilde{g}(\bm{x}).
\end{equation*}
Therefore, by Proposition 8 in \citet{Nolde2022}, we have
\begin{align*}
   \eta^{-1} &= \min_{x_1, x_2 \geq 1}g(x_1,x_2) \\
   &= \min_{x_1, x_2 \geq 1} \min_{x_3,\dots, x_n}\frac{|x_1 - a_{12}x_2 + (a_{12}a_{23}-a_{13})x_3 + \cdots + (a_{12}a_{2n}-a_{1n})x_n|}{1-a_{12}a_{21}} + \\
   &\quad\quad\quad\quad\quad \frac{|x_2 - a_{21}x_1 + (a_{21}a_{13}-a_{23})x_3 + \cdots + (a_{21}a_{1n}-a_{2n})x_n)|} {1-a_{12}a_{21}} + \sum_{i=3}^n |x_i|.
\end{align*}
Using Lemma \ref{optimization} and \ref{optimization12}, we have
\begin{align*}
    \eta &= \Big\{\min_{x_1, x_2 \geq 1} \min_{\substack{i,j=1,\dots,n \\ i\neq j}} \frac{|a_{2i}x_1-a_{1i}x_2|+|a_{2j}x_1-a_{1j}x_2|}{|a_{2i}a_{1j}-a_{1i}a_{2j}|} \Big\}^{-1} \\
    &= \Big\{\min_{\substack{i,j=1,\dots,n \\ i\neq j}} \min_{x_1, x_2 \geq 1}  \frac{|a_{2i}x_1-a_{1i}x_2|+|a_{2j}x_1-a_{1j}x_2|}{|a_{2i}a_{1j}-a_{1i}a_{2j}|} \Big\}^{-1} \\
    &= \Big[\min_{\substack{i,j=1,\dots,n \\ i\neq j}} \min\Big\{ \frac{|a_{2i}-a_{1i}|+|a_{2j}-a_{1j}|}{|a_{2i}a_{1j}-a_{1i}a_{2j}|}, \max(1/a_{1i},1/a_{2i}), \max(1/a_{1j},1/a_{2j}) \Big\} \Big]^{-1}.
\end{align*}
If~\eqref{AI_case}, then by Lemma~\ref{inequalities}, we have $\eta < 1$ and thus $X_1$ and $X_2$ are asymptotically independent.

Notice that although we have assumed $a_{11}=\max_{i\in\{1,\dots,n\}} a_{1i}$ and $a_{22}=\max_{i\in\{1,\dots,n\}} a_{2i}$ at the beginning of this proof, they are only needed to ensure that $a_{11}a_{22}-a_{12}a_{21}\neq 0$, which implies that the matrix $A$ is nonsingular, and they do not affect the expression of $\eta$.
Therefore, the proof is complete.
\end{proof}

%%%%%%%%%%%%%%%%%%%%%%%%%%%%%%%%%%%%% Proofs for Section 4
\vspace*{-10pt}
\section{Proofs for Section \ref{mov_ave_section}} \label{Appendix_proofsfor4}
We first present the proof of Theorem \ref{TypeG_limit_eta_prop}.
\begin{proof}[Proof of Theorem \ref{TypeG_limit_eta_prop}]
By Proposition \ref{AInY}, we know that the residual tail dependence coefficient of $(u_n(\bm{s}_1), u_n(\bm{s}_2))^\top$ is
\[
\eta_n(\bm{s}_1, \bm{s}_2) = \max_{i,j\in \{1,\dots,J\}, i\neq j} \Big\{ \frac{|\Tilde{a}_{1i}\Tilde{a}_{2j} - \Tilde{a}_{1j}\Tilde{a}_{2i}|}{|\Tilde{a}_{1i}-\Tilde{a}_{2i}| + |\Tilde{a}_{1j}-\Tilde{a}_{2j}|}, \Tilde{a}_{1i} \wedge \Tilde{a}_{2i}, \Tilde{a}_{1j}\wedge \Tilde{a}_{2j} \Big\},
%\max_{i,j=1,\dots,J} \Big( &\sup_{\bm{t}_1,\bm{t}_2\in\mathR^2} \frac{|G(\|\bm{t}_1 - \bm{s}_1\|)G(\|\bm{t}_2 - \bm{s}_2\|) - G(\|\bm{t}_1 - \bm{s}_2\|)G(\|\bm{t}_2 - \bm{s}_1\|)|} {G(\bm{0})\{|G(\|\bm{t}_1 - \bm{s}_1\|) - G(\|\bm{t}_1 - \bm{s}_2\|)| + |G(\|\bm{t}_2 - \bm{s}_1\|) - G(\|\bm{t}_2 - \bm{s}_2\|)|\}}, \\
%&\sup_{\bm{t}\in\mathR^2} \frac{G(\|\bm{t} - \bm{s}_1\|)\land G(\|\bm{t} - \bm{s}_2\|)}{G(\bm{0})} \Big)
\]
where $\Tilde{a}_{ri} = \frac{G(\|\bm{s}_r-\bm{d}_i\|)}{G(l_r)},  r=1,2, i=1,\dots,J$ and  $l_r = \min_{i\in\{1,\dots,J\}} \|\bm{s}_r-\bm{d}_i\|$.
If $l_1 \neq l_2$, say $l_1 < l_2$, then one can add one or more mesh nodes close to $\bm{s}_2$ such that $l_1 = l_2^\prime$ holds in the resulting mesh.
So without loss of generality, we assume that $l_1 = l_2 = l$.

For fixed $l>0$, it is clear that $\Tilde{a}_{1i} \wedge \Tilde{a}_{2i} \leq G(\|\bm{s}_1-\bm{s}_2\|/2)/G(l)$, where the equality is achieved when $\bm{d}_i$ is at the midpoint of $\bm{s}_1$ and $\bm{s}_2$, i.e., $\|\bm{s}_1-\bm{d}_i\| = \|\bm{s}_2-\bm{d}_i\| = \|\bm{s}_1-\bm{s}_2\|/2$.
Since our interest is in the case when $l \rightarrow 0$, we assume that $l$ is small and $l < \|\bm{s}_1-\bm{s}_2\|/2$.
We now focus on the term $\frac{|\Tilde{a}_{1i}\Tilde{a}_{2j} - \Tilde{a}_{1j}\Tilde{a}_{2i}|}{|\Tilde{a}_{1i}-\Tilde{a}_{2i}| + |\Tilde{a}_{1j}-\Tilde{a}_{2j}|}$.
Substitute $\Tilde{a}_{ri} = \frac{G(\|\bm{s}_r-\bm{d}_i\|)}{G(l)}$, and denote this term by $f$, i.e.,
$$
f(\bm{d}_1, \bm{d}_2) = \frac{|G(\|\bm{d}_1 - \bm{s}_1\|)G(\|\bm{d}_2 - \bm{s}_2\|) - G(\|\bm{d}_1 - \bm{s}_2\|)G(\|\bm{d}_2 - \bm{s}_1\|)|} {|G(\|\bm{d}_1 - \bm{s}_1\|) - G(\|\bm{d}_1 - \bm{s}_2\|)|G(l) + |G(\|\bm{d}_2 - \bm{s}_1\|) - G(\|\bm{d}_2 - \bm{s}_2\|)|G(l)}.
$$ 
Then we need to find the maximum of this function for $\bm{d}_1, \bm{d}_2 \in \mathR^2\setminus (\mathcal{B}(\bm{s}_1, l) \cup \mathcal{B}(\bm{s}_2, l))$, where $\mathcal{B}(\bm{s}_r, l) = \{\bm{x}\in\mathR^2: \|\bm{x}-\bm{s}_r\| < l\}$ is the open ball of radius $l$ centered at $\bm{s}_r$.

\begin{figure}[!t]
    \centering
    \includegraphics[width=0.5\textwidth]{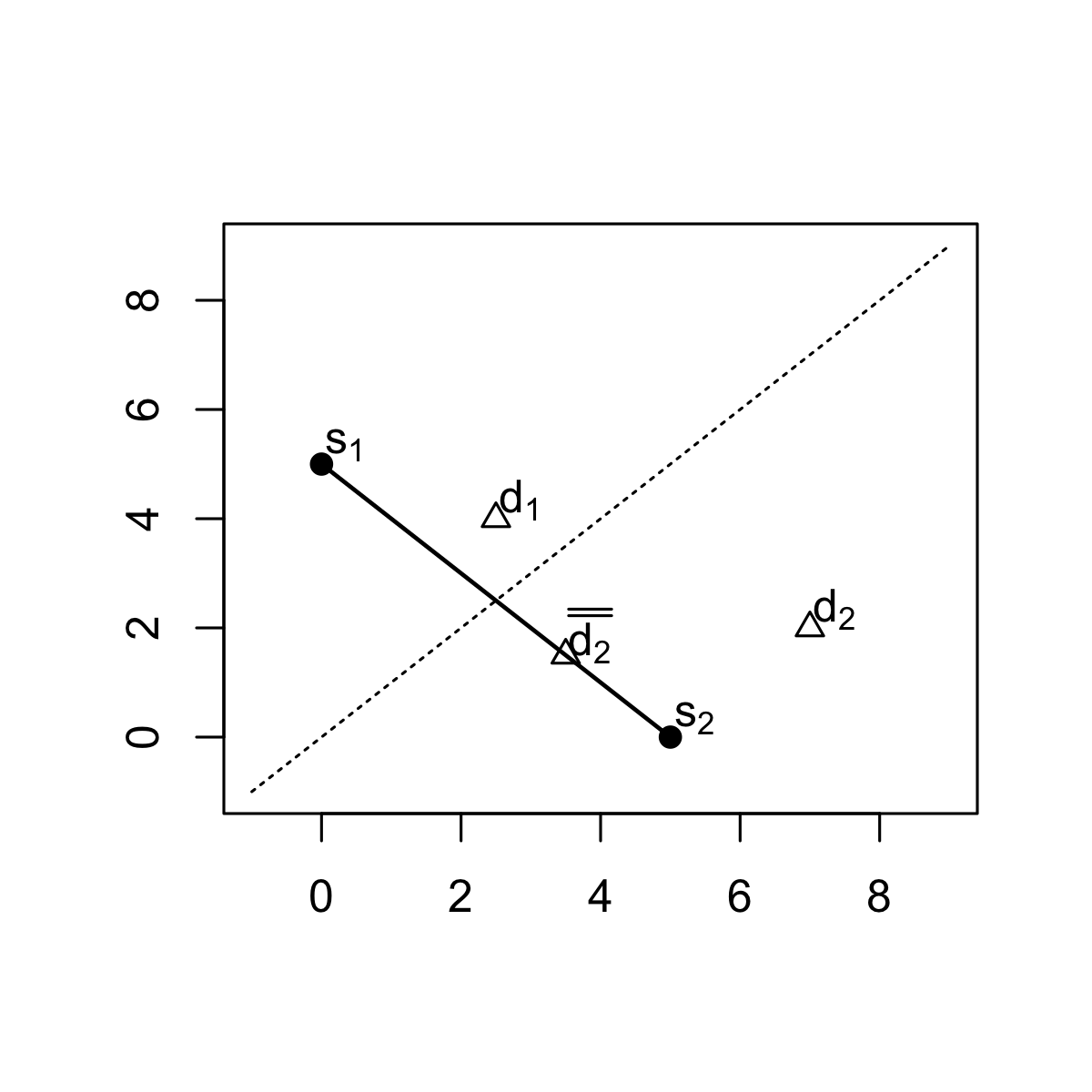}
    \caption{Explaination for the derivation of the limiting $\eta$ in the type G model.}
    \label{figure::typeG_etatrue}
\end{figure}

Denote the line segment from $\bm{s}_1$ to $\bm{s}_2$ by $LS(\bm{s}_1, \bm{s}_2)$, and its perpendicular bisector by $PB(\bm{s}_1, \bm{s}_2)$.
If $\bm{d}_1 \in PB(\bm{s}_1, \bm{s}_2)$ and $\bm{d}_2 \in PB(\bm{s}_1, \bm{s}_2)$, then clearly the function $f$ is not well-defined.
By Proposition \ref{AInY} and the discussion followed, we know that these points can be neglected.
If only one of $\bm{d}_1$ and $\bm{d}_2$ is on the line $PB(\bm{s}_1, \bm{s}_2)$, then we have $f(\bm{d}_1, \bm{d}_2) = G(\|\bm{s}_1 - \bm{s}_2\|/2)/G(l)$.

For $\bm{d}_1, \bm{d}_2 \notin PB(\bm{s}_1, \bm{s}_2)$, without loss of generality we assume that $\bm{d}_1$ is closer to $\bm{s}_1$ than to $\bm{s}_2$, i.e., $\|\bm{d}_1 - \bm{s}_1\| < \|\bm{d}_1 - \bm{s}_2\|$, as shown in Figure \ref{figure::typeG_etatrue}.
One can show that if $\bm{d}_2$ is also closer to $\bm{s}_1$ than to $\bm{s}_2$, then its symmetric point $\Bar{\bm{d}_2}$ about the line $PB(\bm{s}_1, \bm{s}_2)$ satisfies $f(\bm{d}_1, \Bar{\bm{d}_2}) \geq f(\bm{d}_1, \bm{d}_2)$.
This implies that to obtain the maximum of $f$, it is sufficient to consider points $\bm{d}_1, \bm{d}_2$ located on the two sides of the line $PB(\bm{s}_1, \bm{s}_2)$.
Hence we can assume that $\bm{d}_1$ is closer to $\bm{s}_1$ than to $\bm{s}_2$, and $\bm{d}_2$ is closer to $\bm{s}_2$ than to $\bm{s}_1$.
Since $G$ is strictly decreasing on $[0,\infty)$, the function $f$ can then be simplified to
\[
f(\bm{d}_1, \bm{d}_2) = \frac{G(\|\bm{d}_1 - \bm{s}_1\|)G(\|\bm{d}_2 - \bm{s}_2\|) - G(\|\bm{d}_1 - \bm{s}_2\|)G(\|\bm{d}_2 - \bm{s}_1\|)} {\{G(\|\bm{d}_1 - \bm{s}_1\|) - G(\|\bm{d}_1 - \bm{s}_2\|) + G(\|\bm{d}_2 - \bm{s}_1\|) - G(\|\bm{d}_2 - \bm{s}_2\|)\}G(l) }.
\]

Furthermore, due to the convexity of $G$, one can show that for $\bm{d} \in \mathR^2\setminus (\mathcal{B}(\bm{s}_1, l) \cup \mathcal{B}(\bm{s}_2, l))$, the range of the function $|G(\|\bm{d} - \bm{s}_1\|) - G(\|\bm{d} - \bm{s}_2\|)|$ is $[0, G(l)-G(\|\bm{s}_1 - \bm{s}_2\| - l)]$ with its minimum obtained on the line $PB(\bm{s}_1, \bm{s}_2)$ and its maximum at the intersection points of the line passing through $\bm{s}_1$ and $\bm{s}_2$ and the two circles with radius $l$ centered at $\bm{s}_1$ and  $\bm{s}_2$.
Now if we fix the point $\bm{d}_1$, then by the intermediate value theorem \citep[Theorem 9.7.1]{Tao2016Analysis}, for any $\bm{d}_2$ such that ${\|\bm{d}_2 - \bm{s}_2\| < \|\bm{d}_2 - \bm{s}_1\|}$ and ${\bm{d}_2\notin \mathcal{B}(\bm{s}_2, l)}$, there always exists one point $\Bar{\Bar{\bm{d}_2}}$ on the line segment ${LS((\bm{s}_1+\bm{s}_2)/2, \bm{s}_2)\setminus \mathcal{B}(\bm{s}_2, l)}$ such that ${G(\|\bm{d}_2 - \bm{s}_1\|) - G(\|\bm{d}_2 - \bm{s}_2\|) =  G(\|\Bar{\Bar{\bm{d}_2}} - \bm{s}_1\|) - G(\|\Bar{\Bar{\bm{d}_2}} - \bm{s}_2\|)}$, and ${G(\|\Bar{\Bar{\bm{d}_2}} - \bm{s}_i\|) \geq G(\|\bm{d}_2 - \bm{s}_i\|)}$, $i=1,2$, which yields ${f(\bm{d}_1, \Bar{\Bar{\bm{d}_2}}) \geq f(\bm{d}_1, \bm{d}_2)}$.
Hence, it is sufficient to consider only points $\bm{d}_1 \in LS(\bm{s}_1, (\bm{s}_1+\bm{s}_2)/2)\setminus\mathcal{B}(\bm{s}_1, l)$ and $\bm{d}_2 \in LS(\bm{s}_2, (\bm{s}_1+\bm{s}_2)/2)\setminus\mathcal{B}(\bm{s}_2, l)$, and the optimization problem can be rewritten as
\[
\sup_{x_1, x_2 \in [l,h)} \Bar{f}(x_1,x_2) =  \sup_{x_1, x_2 \in [l,h)} \frac{G(x_1)G(x_2) - G(2h-x_1)G(2h-x_2)}{\{G(x_1)-G(2h-x_1)+G(x_2)-G(2h-x_2)\}G(l)}.
\]

Since $G$ is absolutely continuous, we have the partial derivative of $\Bar{f}$ with respect to $x_1$ as
\[
\frac{\partial \Bar{f}}{\partial x_1} = a(x_1, x_2) [G^\prime(x_1)\{G(x_2)-G(2h-x_2)\} - G^\prime(2h-x_1)\{G(x_1)-G(2h-x_2)\} ],
\]
where $a(x_1,x_2)=\frac{G(x_2)-G(2h-x_2)}{\{G(x_1)-G(2h-x_1)+G(x_2)-G(2h-x_2)\}^2 G(l)} > 0$.
For $x\in [l,h)$, due to the convexity of $G$, we have that $G^\prime(x) \leq G^\prime(2h-x)$.
Hence, 
$$\frac{\partial \Bar{f}}{\partial x_1} \leq a(x_1,x_2)G^\prime(2h-x_1)\{G(x_2)+G(2h-x_2)-G(x_1)-G(2h-x_1)\}. $$
As $G$ is monotonically decreasing, we know that $G^\prime(2h-x_1)\leq 0$.
Furthermore, the convexity of $G$ implies that the function $g(x)=G(x)+G(2h-x)$ is monotonically decreasing on $[l,h)$.
Therefore, for fixed $x_2\in [l, h)$, we know that if $x_1\geq x_2$, then $\frac{\partial \Bar{f}}{\partial x_1} \leq 0$.
This implies that 
\[
\Bar{f}(x_1,x_2) \leq \Bar{f}(x_2,x_2) = \{G(x_2)+G(2h-x_2)\}/\{2G(l)\} \leq \{G(l)+G(2h-l)\}/\{2G(l)\},
\]
where the maximum is obtained when $x_1=x_2=l$.

Above all, we have that for $\bm{d}_1, \bm{d}_2 \in \mathR^2\setminus (\mathcal{B}(\bm{s}_1, l) \cup \mathcal{B}(\bm{s}_2, l))$,
\begin{align*}
    f(\bm{d}_1, \bm{d}_2) \leq \max\Big\{\frac{G(\|\bm{s}_1 - \bm{s}_2\|/2)}{G(l)}, \frac{G(l)+G(\|\bm{s}_1 - \bm{s}_2\|-l)}{2G(l)} \Big\} = 1/2 + G(\|\bm{s}_1 - \bm{s}_2\|-l)/G(l),
\end{align*}
where the last equality holds due to the convexity of $G$, and the maximum of $f$ is achieved when $\bm{d}_1, \bm{d}_2$ are the intersection points of the line segment $LS(\bm{s}_1,\bm{s}_2)$ and the two circles with radius $l$ centered at $\bm{s}_1$ and $\bm{s}_2$.
Hence, as the mesh becomes finer and the set of mesh nodes $M_n$ becomes denser in $\mathcal{D}$, which means that $l\rightarrow 0$, the residual tail dependence function $\eta(\bm{s}_1,\bm{s}_2)$ tends to $1/2 + G(\|\bm{s}_1 - \bm{s}_2\|)/G(0)$.
This completes the proof.
\end{proof}

We now prove Theorem~\ref{OU_limit_eta}.
\begin{proof}[Proof of Theorem~\ref{OU_limit_eta}]
By Proposition~\ref{AInY}, the residual tail dependence coefficient of the approximation $(u_n(s_1), u_n(s_2))$ of the form~\eqref{mov_OU_approx} is
\[
\eta = \max_{i,j=0,\dots,n_2-1, i\neq j} \Big\{\frac{|a_{2i}a_{1j}-a_{1i}a_{2j}|}{|a_{2i}-a_{1i}|+|a_{2j}-a_{1j}|}, \min(a_{1i},a_{2i}), \min(a_{1j},a_{2j}) \Big\},
\]
with ${a_{1i}=G(s_1-t_i)/G(s_1-t_{n_1-1})}$ for $0\leq i\leq n_1-1$ and ${a_{1i}=0}$ for ${n_1\leq i\leq n_2-1}$, and ${a_{2i}=G(s_2-t_i)/G(s_2-t_{n_2-1})}$ for ${0\leq i\leq n_2-1}$.

Since $G$ is strictly decreasing, the sequence ${\{a_{1i},0\leq i\leq n_1-1\}}$ and ${\{a_{2i},0\leq i\leq n_2-1\}}$ are strictly increasing.
Hence, ${\max_{i=0,\dots,n_2-1}\min(a_{1i},a_{2i}) = a_{2(n_1-1)}}$.
For the first term inside the maximum operation of the expression of $\eta$, we have 
\begin{align*}
    \frac{|a_{2i}a_{1j}-a_{1i}a_{2j}|}{|a_{2i}-a_{1i}|+|a_{2j}-a_{1j}|} &\leq \frac{\max(a_{1j}, a_{2j}) |a_{2i}-a_{1i}|}{|a_{2i}-a_{1i}|+|a_{2j}-a_{1j}|} 
    = \frac{\max(a_{1j}, a_{2j})}{1+|a_{2j}-a_{1j}|/|a_{2i}-a_{1i}|} \\
    &\leq \frac{\max(a_{1j}, a_{2j})}{1+|a_{2j}-a_{1j}|}
    = \frac{1}{1+\{1-\min(a_{1j},a_{2j})\}/\max(a_{1j}, a_{2j})} \\
    &\leq \frac{1}{2-a_{2(n_1-1)}},
\end{align*}
where the equality is obtained when $i=n_2-1$ and $j=n_1-1$.
Thus we have $\eta$ for $(u_n(s_1), u_n(s_2))^\top$ as 
\[
\eta = \max\Big(\frac{1}{2-a_{2(n_1-1)}}, a_{2(n_1-1)} \Big) = \frac{1}{2-a_{2(n_1-1)}}.
\]

Therefore, as $m_n\rightarrow\infty$ and the set of mesh nodes $M_n$ tries to fill up the space in $\mathcal{D}$, we have the limiting residual tail dependence function as $1/\{2-G(h)/G(0)\}$.
\end{proof}

Next, we prove Theorem \ref{OUprocess_AI}.
\begin{proof}[Proof of Theorem \ref{OUprocess_AI}]
Without loss of generality we assume that $0 \leq s_1 < s_2$.
The representation (\ref{OUrepresentation}) implies that $u(s_2)$ can be represented as the convolution of two independent random variables
\[
  u(s_2) = e^{-a(s_2-s_1)}u(s_1) + \int_{s_1}^{s_2} e^{-a(s_2-s)}{\mathrm d}z(at).
\]
Since $u(s) = V$ and $\bar{F}_V \in \mathscr{L}_\beta$, we know that $u(s_1), u(s_2)$ have exponential tails with index $\beta$ and $e^{-a(s_2-s_1)}u(s_1)$ has an exponential tail with index $e^{a(s_2-s_1)}\beta$.
Denote $V_1=e^{-a(s_2-s_1)}u(s_1)$.
The self-decomposability of $V_1$ yields its infinite divisibility, which further implies that the Wiener condition is satisfied, i.e. $M_{V_1}(\beta+it)\neq 0, t\in\mathR$, where $M_{V_1}$ is the moment generating function of $V_1$; see Theorem 25.17 of \citep{Sato1999} for more details.
Then Lemma 2.5 in \citet{Pakes2004} yields that $\int_{s_1}^{s_2} e^{-a(s_2-s)}{\mathrm d}z(at)$ must have an exponential tail with index $\beta$.
Therefore, the result in Proposition \ref{AInY} and
Example \ref{AItwoY} gives the asymptotic independence between $u(s_1)$ and $u(s_2)$, and their residual tail dependence coefficient is $\eta = 1/(2-e^{-a(s_2-s_1)})$.    
\end{proof}

In the following we prove Proposition \ref{OUexpotail}.
\begin{proof}[Proof of Proposition \ref{OUexpotail}]
By \citet{Sgibnev1990} and \citet{Shimura2005}, we know that an infinitely divisible distribution is convolution tail equivalent if and only if its normalized L\'evy measure is convolution tail equivalent.
That is, for any $t\in\mathR$,
    \[
    \Bar{F}_{z(1)} \in\mathscr{S}_\beta
    \quad \quad \Leftrightarrow \quad \quad
    \frac{\mathds{1}_{x>1} U_{z(1)}([x,\infty))}{U_{z(1)}([1,\infty))} \in\mathscr{S}_\beta.
    \]
From \citet{BarndorffNielsen1998} we know that if $U_{z(1)}([x,\infty))$ is continuous in $x$ on $(0,\infty))$, then we have the relation $U_{u}([x,\infty)) = \int_x^\infty s^{-1} U_{z(1)}([s,\infty)) {\mathrm d}s$.
Hence, for $x>1$,
    \[
    U_{u}([\log x,\infty)) = \int_{\log x}^\infty s^{-1} U_{z(1)}([s,\infty)) {\mathrm d}s = \int_x^\infty U_{z(1)}([\log t,\infty))/(t\log t) {\mathrm d}t.
    \]
Since $\frac{\mathds{1}_{x>1} U_{z(1)}([x,\infty))}{U_{z(1)}([1,\infty))} \in\mathscr{S}_\beta$, we have that $U_{z(1)}([\log t,\infty)) \in{\text RV}_{-\beta}$, which further implies that ${U_{z(1)}([\log t,\infty))/(t\log t) \in{\text RV}_{-\beta-1}}$.
By Karamata's theorem (cf. Theorem 2.1 in \citet{Resnick2007}), 
    \[
    \lim_{x\rightarrow\infty} \frac{U_{z(1)}([\log x,\infty))/\log x}{\int_x^\infty U_{z(1)}([\log t,\infty))/(t\log t) {\mathrm d}t} = \beta.
    \]
Hence, 
    \[
    \lim_{x\rightarrow\infty} \frac{U_{u}([\log x,\infty))}{U_{z(1)}([\log t,\infty))} = \lim_{x\rightarrow\infty} \frac{1}{\beta \log x} = 1/\beta.
    \]
Therefore, using Theorem 1 in \citet{Cline1986}, we have that $\frac{\mathds{1}_{x>1} U_{u}([x,\infty))}{U_{u}([1,\infty))} \in\mathscr{S}_\beta$ and thus $\Bar{F}_{u}\in\mathscr{S}_\beta$.

Now we denote the L\'evy measure of $z(t)$ as $U_{z(t)}$.
Since $F_{z(1)}\in \mathscr{S}_\beta$, we know that ${\frac{\mathds{1}_{x>0} U_{z(1)}}{U_{z(1)}([1,\infty))} \in\mathscr{S}_\beta}$.
By the definition of L\'evy process, we have $U_{z(t)}=t U_{z(1)}$, for any $t>0$.
Hence, using Theorem 1 in \citet{Cline1986}, we know that $\frac{\mathds{1}_{x>0} U_{z(t)}([x,\infty))}{U_{z(t)}([1,\infty))} \in\mathscr{S}_\beta$ and thus $\Bar{F}_{z(t)}\in \mathscr{S}_\beta$.
\end{proof}

%%%%%%%%%%%%%%%%%%%%%%%%%%%%%%%%%%%%% FEM approximation of SPDEs
\vspace*{-10pt}
\section{Finite Element Approximations}\label{Appendix_FEM}
In this section we introduce the finite element approximation which is the key for efficient simulation and inference of the type G Mat\'ern fields \citep{Bolin2014,BolinWallin2020}.

Similarly to the non-Gaussian OU processes, the type G Mat\'ern SPDE random fields are an important extension of the well-known SPDE-based formulation of Gaussian random fields \citep{Lindgren2011}, aiming to capture more flexible marginal behaviors, different sample path properties, and non-Gaussian dependence structures.
Specifically, the generalization from Gaussian white noise to general L\'evy noise in the SPDE (\ref{TypeGmodel}) provides a natural way to construct non-Gaussian random fields, while the differential operator $(\kappa^2 - \Delta) ^{\alpha/2}$ ensures the Mat\'ern covariance for the resulting processes.

The main advantage of the SPDE-based representation is that one can approximate its solution by using the finite element method to achieve computationally efficient simulation and inference.
To make sense of the SPDE (\ref{TypeGmodel}), one can think of $\mathcal{M}$ as a homogeneous L\'evy basis, namely an infinitely divisible and independently scattered random measure, and interpret the Equation (\ref{TypeGmodel}) in a weak sense.
That is, for a test function $g$ in an appropriate space, we have
\[
\langle g(\bm{s}), (\kappa^2 - \Delta) ^{\alpha/2} u(\bm{s}) \rangle = \int g(\bm{s}) \mathcal{M}({\mathrm d}\bm{s}),
\]
where the equality holds in distribution, and $\langle a(\bm{s}), b(\bm{s}) \rangle = \int a(\bm{s}) b(\bm{s}) {\mathrm d}\bm{s}$.
To approximate the weak solution numerically efficiently, we consider a bounded domain $\mathcal{D}\in\mathR^d$ and assume that the operator $(\kappa^2 - \Delta) ^{\alpha/2}$ is equipped with suitable boundary conditions \citep{Lindgren2011,Bolin2014}.
Then one can choose the space of the test functions to be the linear span of a finite number of basis functions $\{\phi_i, i=1,\dots,n\}$.
For instance, similarly as in \citet{Lindgren2011}, $\phi_i$ can be chosen as the commonly used piecewise linear basis functions induced by a triangulation of the domain $\mathcal{D}$.
We then obtain a finite element approximation of the weak solution as
\begin{equation*}
    u_n(\bm{s}) = \sum_{i=1}^n \omega_i \phi_i(\bm{s}),
\end{equation*}
where the stochastic weights $\bm{\omega}=(\omega_1,\dots,\omega_n)^\top$ satisfy
\[
    K_\alpha \bm{\omega} = \Big(\int_\mathcal{D} \phi_1(\bm{s}) \mathcal{M}({\mathrm d}\bm{s}), \dots, \int_\mathcal{D} \phi_n(\bm{s}) \mathcal{M}({\mathrm d}\bm{s}) \Big)^\top,
\]
with $K_\alpha = C\{C^{-1}(\kappa^2 C + G)\}^{\alpha/2}$, and the matrices $C, G \in\mathR^{n\times n}$ having elements $C_{ij}=\langle\phi_i, \phi_j\rangle_\mathcal{D}$, and $G_{ij}=\langle\nabla\phi_i, \nabla\phi_j\rangle_\mathcal{D}$, respectively.

The random measure $\mathcal{M}$ is said to be of type G if $\mathcal{M}$ evaluated on the unit square $[0,1]^d$ can be represented as a normal mean-variance mixture, i.e., $\mathcal{M}([0,1]^d) = \mu + \gamma v + \sqrt{v}Z$, where $v$ is a non-negative infinitely divisible random variable and $Z$ is a Gaussian random variable.
This representation of the type G random measure admits a Gaussian distribution for $\bm{\omega}$ conditioning on the variance process $v(\bm{s})$. 
By assuming that $v$ is closed under convolution, the distribution of $\bm{\omega}$ can be approximated by
\begin{equation*}
    \bm{\omega} \mid \Bar{\bm{v}} \sim \mathcal{N}\big(K_\alpha^{-1} (\mu + \gamma \Bar{\bm{v}}), K_\alpha^{-1} \text{diag}(\Bar{\bm{v}}) K_\alpha^{-1}\big)
\end{equation*}
with ${\bar{\bm{v}}=(\bar{v}_1,\dots,\bar{v}_n)^\top}$ and ${\bar{v}_j = M_v(\mathcal{D}_j)}$, where $M_v$ is the random measure associated with $v$, and ${\mathcal{D}_j = \{\bm{s}: \phi_j(\bm{s}) \geq \phi_i(\bm{s}), \forall i\neq j \}}$.
For the GIG distribution introduced in Section \ref{ExpoIntro}, the only two subclasses which are closed under convolution are the Gamma distribution and the inverse Gaussian distribution, which lead to the variance Gamma and normal inverse Gaussian random measures, respectively, for $\mathcal{M}$, and these are the two non-Gaussian SPDE models considered in \citet{Bolin2014}, \citet{WallinBolin2015} and \citet{BolinWallin2020}.

If $v$ has an inverse Gaussian distribution, then $\bar{v}_i, i=1,\dots,n$, are independent and inverse Gaussian distributed, and have exponential tails with the same index.
Consequently, the approximate weak solution $u_n(\bm{s})$ can be rewritten as
\begin{align}
    u_n(\bm{s}) = (\phi_1(\bm{s}), \dots, \phi_n(\bm{s})) \bm{\omega}
    = (\phi_1(\bm{s}), \dots, \phi_n(\bm{s})) K_\alpha^{-1} \bm{Y}, \label{FEmapprox}
\end{align}
where $\bm{Y}=(Y_1,\dots,Y_n)^\top$, $Y_i,i=1,\dots,n$ are independent and normal inverse Gaussian distributed, and have exponential tails with the same index.
If $v$ has a Gamma distribution, then we have the same representation as in (\ref{FEmapprox}) for the approximate weak solution, but with the difference that $Y_i,i=1,\dots,n$ have a variance Gamma distribution.

The preceding analysis shows that characterizing the extremal dependence of the finite element approximation $u_n(\bm{s})$ is equivalent to characterizing the extremal dependence of the random vector $\bm{X}=(X_1, X_2)^\top$ constructed as in (\ref{discretemodel}).
From the results in Section \ref{dissection} we know that the extremal dependence between $u_n(\bm{s}_1)$ and $u_n(\bm{s}_2)$ depends on whether $\argmax_{i=1,\dots,n} \{(\phi_1(\bm{s}_1), \dots, \phi_n(\bm{s}_1)) K_\alpha^{-1}\}$ is equal to $\argmax_{i=1,\dots,n} \{(\phi_1(\bm{s}_2), \dots, \phi_n(\bm{s}_2)) K_\alpha^{-1}\}$.
And for any given mesh and basis functions $\{\phi_i\}$, the residual tail dependence coefficient of the finite element approximation can be easily computed using the formula in Proposition~\ref{AInY}.

\vspace*{-10pt}
\section{Additional Simulation Studies}\label{Appendix_simulation}

\begin{figure}[!t]
    \centering
    \includegraphics[width=0.8\textwidth]{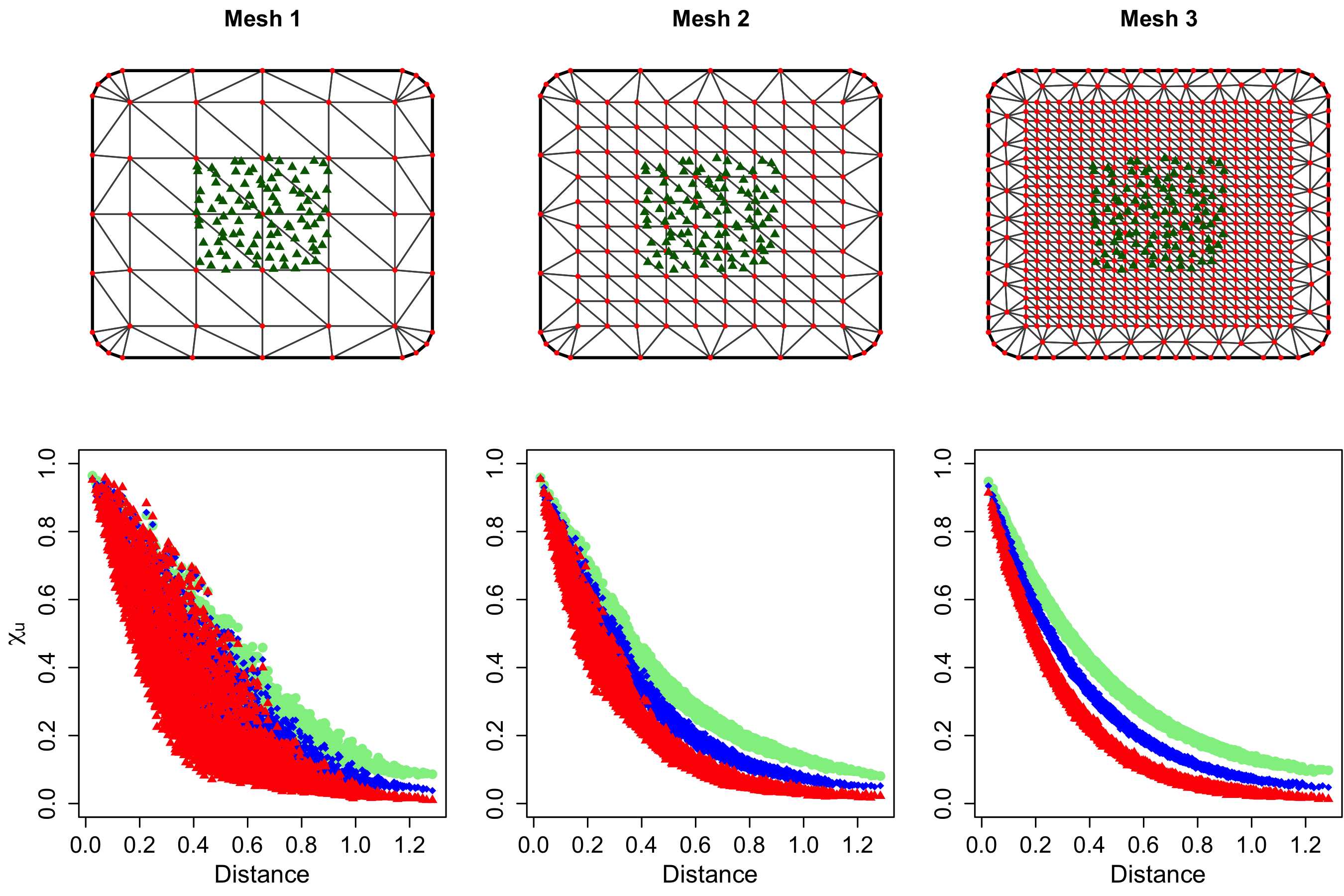}
    \caption{Extremal dependence of the finite element approximations of the NIG Mat\'ern SPDE model with respect to different mesh constructions. Mesh $1$, $2$, and $3$ are constructed on a lattice with respective $25$, $100$, and $625$ nodes. The empirical tail dependence coefficient $\chi(q)$ with different $q$ (green points: $q=0.95$; blue points: $q=0.975$; red points: $q=0.99$) of the finite element approximation model based on mesh $1$, $2$, and $3$ are shown in the lower left, middle, and right panels, respectively.}
    \label{figure::mesh_illustration}
\end{figure}
Here we examine the effect of different mesh constructions on the extremal dependence of the induced finite element approximation of NIG Mat\'ern SPDE models.
We choose the same SPDE parameters and NIG noise parameters as in the first simulation study in Section~\ref{TypeG_section}.
As shown in Figure \ref{figure::mesh_illustration}, we randomly select 100 sites in the unit square and consider three different mesh constructions covering the sites.
The meshes are constructed based on lattices with outer extensions, where meshes 1, 2, and 3 have 25, 100, and 625 lattice nodes, respectively.
We then simulate $10^5$ observations from the finite element approximation model and the lower panel in Figure \ref{figure::mesh_illustration} depicts the empirical tail dependence coefficient of all pairs of sites with respect to different quantile levels $u=0.95$ (green points), $0.975$ (blue points), and $0.99$ (red points).
This figure clearly illustrates the nonstationarity induced by a coarse mesh (mesh 1 and 2), i.e., various values of $\chi(q)$ exist for the same distance, but this nonstationarity tends to vanish as the mesh becomes finer.
One can also observe in the lower right panel that as the quantile level $u$ increases to $1$, $\chi(q)$ tends to decrease.
This seems to coincide with our theoretical findings of the integral approximation in Section \ref{TypeG_section}, whose limiting extremal dependence class is asymptotic independence and which necessarily implies that $\chi=\lim_{q\uparrow 1}\chi(q)=0$ for all pairs.

%%%%%%%%%%%%%%%%%%%%%%%%%%%%%%%%%%%%% Proofs for auxiliary results
\vspace*{-10pt}
\section{Additional Results}\label{Appendix_aux_results}

%\noindent\textbf{Proposition C.1.}
\begin{proposition}\label{Aux_result_1}
Let $\chi$ be the tail dependence coefficient in Example \ref{ADtwoY}. If $a_{12}\neq a_{22}$, then $\chi < 1$.
\end{proposition}
\begin{proof}
Without loss of generality, we assume that $0\leq a_{12} < a_{22} <1$.
To prove $\chi < 1$, we first show that $c > 0$.
Note that for $0 < t < \sqrt{\psi}$, we have
\begin{align*}
M^\prime_{Y_1}(t) &= \int_{-\infty}^\infty y\exp(ty) F_{Y_1}({\mathrm d}y) \\
&= \int_{-\infty}^0 y\exp(ty) F_{Y_1}({\mathrm d}y) + \int_{0}^\infty y\exp(ty) F_{Y_1}({\mathrm d}y) \\
&>  \int_{-\infty}^0 y  F_{Y_1}({\mathrm d}y) + \int_{0}^\infty y  F_{Y_1}({\mathrm d}y)
= 0,
\end{align*}
where the last equality is due to the symmetry of the distribution ${\text{GH}}(\lambda, \tau, \psi, \mu=0, \gamma=0)$.
Hence, $c=\{\log M_{Y_1}(a_{22}\sqrt{\psi}) - \log M_{Y_1}(a_{12}\sqrt{\psi}) \}/\sqrt{\psi} > 0$.

Therefore,
\begin{align*}
    1 - \chi &= 1-\frac{\int_{(a_{22}-a_{12})y\leq c} \exp(a_{22}\sqrt{\psi}y)  F_{Y_1}({\mathrm d}y)}{M_{Y_1}(a_{22}\sqrt{\psi})} - 
    \frac{\int_{(a_{22}-a_{12})y> c} \exp(a_{12}\sqrt{\psi}y)  F_{Y_1}({\mathrm d}y)}{M_{Y_1}(a_{12}\sqrt{\psi})} \\
    &= \int_{c/(a_{22}-a_{12})}^{\infty} e^{\sqrt{\psi}a_{22}y - \log M_{Y_1}(a_{22}\sqrt{\psi})} - e^{\sqrt{\psi}a_{12}y - \log M_{Y_1}(a_{12}\sqrt{\psi})} F_{Y_1}({\mathrm d}y) \\
    &>0,
\end{align*}
where the last inequality holds since ${\sqrt{\psi}a_{22}y - \log M_{Y_1}(a_{22}\sqrt{\psi}) > \sqrt{\psi}a_{12}y - \log M_{Y_1}(a_{12}\sqrt{\psi})}$ for ${y > c/(a_{22}-a_{12})}$.
\end{proof}

%\noindent\textbf{Proposition C.2.}
%\begin{proposition}\label{Aux_result_2}
%    Let $z(t)$ be a L\'evy process. Denote $F_{z(t)}$ as the distribution function of $z(t)$. If $z(1)$ has a convolution tail equivalent distribution, i.e. $F_{z(1)}\in \mathscr{S}_\beta, \beta\geq 0$, then we have $F_{z(t)}\in \mathscr{S}_\beta, \beta\geq 0$.
%\end{proposition}
%\begin{proof}
%Denote the L\'evy measure of $z(t)$ as $U_{z(t)}$.
%Since $F_{z(1)}\in \mathscr{S}_\beta$, by \citet{Sgibnev1990} and \citet{Shimura2005}, we know that $\frac{\mathds{1}_{x>0} U_{z(1)}}{U_{z(1)}([1,\infty))} \in\mathscr{S}_\beta$.
%By the definition of L\'evy process, we have $U_{z(t)}=t U_{z(1)}$, for any $t>0$.
%Hence, using Theorem 1 in \citet{Cline1986}, we know that $\frac{\mathds{1}_{x>0} U_{z(t)}}{U_{z(t)}([1,\infty))} \in\mathscr{S}_\beta$ and thus $F_{z(t)}\in \mathscr{S}_\beta$.
%\end{proof}

%\noindent\textbf{Proposition C.3.}
\begin{proposition}\label{Aux_result_3}
    Let $f(x) = x^\nu K_\nu(x), \nu>0, x>0$. Then $f(x)$ is strictly decreasing on $(0,\infty)$.
\end{proposition}
\begin{proof}
Using the formula $K_\nu^\prime(x)=-\{K_{\nu-1}(x)+K_{\nu+1}(x)\}/2$ (cf. \citet{Abramowitz1972}, formula 9.6.26) we have
\begin{align*}
    f^\prime(x) &= \nu x^{\nu-1}K_\nu(x) + x^\nu \Big\{-\frac{K_{\nu-1}(x)+K_{\nu+1}(x)}{2} \Big\} \\
    &= x^{\nu-1} \frac{2\nu K_\nu(x)-xK_{\nu-1}(x)-xK_{\nu+1}(x)}{2}.
\end{align*}
Since $xK_{\nu+1}(x)=xK_{\nu-1}(x)+2\nu K_\nu(x)$ (cf. \citet{Abramowitz1972}, formula 9.6.26), we have 
\begin{align*}
    f^\prime(x) &= x^{\nu-1} \frac{2\nu K_\nu(x)-xK_{\nu-1}(x)-xK_{\nu-1}(x)-2\nu K_\nu(x)}{2}
    = -x^\nu K_{\nu-1}(x)
    < 0.
\end{align*}
Therefore, $f(x)$ is strictly decreasing on $(0,\infty)$.
\end{proof}

We now give more details of the GH distribution.
The Lebesgue density function of ${\text{GH}}(\lambda,\tau,\psi,\mu,\gamma)$ is
\[
    f_{\text{GH}}(x) = \frac{\psi^{\lambda/2} (\psi+\gamma^2)^{1/2 - \lambda}}{(2\pi)^{1/2}\tau^{\lambda/2}K_\lambda(\sqrt{\tau \psi})} \cdot \frac{K_{\lambda-1/2}\Big[\sqrt{\{\tau+(x-\mu)^2\}(\psi+\gamma^2)}\Big] e^{\gamma(x-\mu)}} {\Big[\sqrt{\{\tau+(x-\mu)^2\}(\psi+\gamma^2)} \Big]^{1/2 - \lambda}}, \quad x\in\mathR.
\]
To see why the GH density function also has exponential tails, using the asymptotic expansion of the modified Bessel function of the second kind for large arguments \citep[Formula 9.7.2]{Abramowitz1972}, we obtain
\begin{equation*}
    f_{\text{GH}}(x) \sim c_0 x^{\lambda-1}e^{-(\sqrt{\psi+\gamma^2}-\gamma)x}, \quad x\rightarrow \infty,
\end{equation*}
where $c_0>0$ is a constant, and $f(x)\sim g(x)$ as $x\rightarrow \infty$ means that $g$ is eventually non-zero and $f(x)/g(x) \rightarrow 1$ as $x\rightarrow \infty$.
A simple application of L'H\^opital's rule gives that the GH distribution functions also have exponential tails.

\end{appendices}

\bibliography{reference}

\end{document}